\documentclass[journal]{IEEEtran}
\usepackage{amsmath,amsfonts}
\usepackage{algorithmic}
\usepackage{algorithm}
\usepackage{array}
\usepackage[caption=false,font=normalsize,labelfont=sf,textfont=sf]{subfig}
\usepackage{textcomp}
\usepackage{stfloats}
\usepackage{url}
\usepackage{verbatim}
\usepackage{cite}
\usepackage{fullpage}
\usepackage[T1]{fontenc}
\usepackage{color}
\usepackage{babel}
\usepackage{verbatim}
\usepackage{float}
\usepackage{amsmath}
\usepackage{amsthm}
\usepackage{graphicx}
\usepackage{mathrsfs}
\usepackage{amssymb}
\usepackage{wasysym}
\usepackage{microtype}
\usepackage{silence}
\usepackage{xcolor}

\theoremstyle{plain}
\newtheorem{thm}{Theorem}[section]
\newtheorem{lem}{Lemma}[section]

\newtheorem{assum}{Assumption}[section]

\makeatother
\DeclareMathOperator{\Expect}{\mathbb{E}}

\graphicspath{{figures/}}
\begin{document}

\title{Convergence Analysis of the Last Iterate in Distributed Stochastic Gradient Descent with Momentum}

\author{Difei~Cheng,~\IEEEmembership{Member,~IEEE,} Ruinan~Jin,~\IEEEmembership{Member,~IEEE,}  Hong~Qiao,~\IEEEmembership{Fellow,~IEEE,} and Bo~Zhang,~\IEEEmembership{Senior Member,~IEEE}
\thanks{D. Cheng and H. Qiao is with State Key Laboratory of Multimodal Artificial Intelligence Systems,
	Institute of Automation, Chinese Academy of Sciences, Beijing 100190, China (email:chengdifei@amss.ac.cn~, hong.qiao@ia.ac.cn)}  
\thanks{R. Jin is with School of Data Science, The Chinese University of Hong Kong at Shenzhen, Shenzhen 518172, China (email: jinruinan@cuhk.edu.cn)}       
\thanks{B. Zhang is with LSEC and Institute of Applied Mathematics, AMSS,
	Chinese Academy of Sciences, Beijing 100190, China (email: b.zhang@amt.ac.cn)}
\thanks{Corresponding author: Ruinan Jin}
}

\markboth{Journal of \LaTeX\ Class Files,~Vol.~14, No.~8, August~2021}%
{Shell \MakeLowercase{\textit{et al.}}: A Sample Article Using IEEEtran.cls for IEEE Journals}

\maketitle

\begin{abstract}

Distributed stochastic gradient methods are widely used to preserve data privacy and ensure scalability in large-scale learning tasks. While existing theory on distributed momentum Stochastic Gradient Descent (mSGD) mainly focuses on time-averaged convergence, the more practical last-iterate convergence remains underexplored. In this work, we analyze the last-iterate convergence behavior of distributed mSGD in non-convex settings under the classical Robbins–Monro step-size schedule. We prove both almost sure convergence and  $L_2$ convergence of the last iterate, and derive convergence rates. We further show that momentum can accelerate early-stage convergence, and provide experiments to support our theory.
\end{abstract}

\section{Introduction}
Stochastic gradient  methods ~\cite{1951A,1964Some,tang2018adaptive,kim2014combining} play a central role in machine learning due to their efficiency in minimizing regularized empirical risk, including in the training of deep neural networks~\cite{Graves2013Speech,Nguyen2018SGD,HintonReducing}. 

Classical stochastic gradient methods are typically designed for \textbf{centralized architectures}, where a central server collects data from multiple edge devices (workers) to compute gradients. However, such centralized training frameworks face several challenges, including \textbf{privacy concerns} arising from raw data sharing and \textbf{limited scalability} due to computational bottlenecks.

To overcome these limitations, \textbf{distributed stochastic gradient methods} have been developed. In such settings, each worker computes stochastic gradients using local data, and model parameters are periodically synchronized across nodes. This paradigm solves optimization problems of the form:
\begin{align}\label{distributed_sgd}
	\min_{x\in \mathbb{R}^{N}} \mathbb{E}_{\xi} \left[g(x,\xi)\right], \quad \text{with} \quad g(x) = \frac{1}{m} \sum_{i=1}^{m} g_i(x),
\end{align}
where $g(x, \xi)$ is an unbiased stochastic estimator of the objective function $g(x)$, $\xi$ denotes stochastic sampling noise or perturbations, and $m$ is the number of worker nodes.

Early distributed algorithms typically rely on frequent model synchronization, which incurs substantial communication overhead due to bandwidth limitations, network latency, or hardware constraints—challenges that become particularly severe when training deep models. To alleviate such communication costs and enhance algorithmic efficiency, several classical variants of distributed stochastic gradient algorithms have been proposed:

\begin{itemize}
    \item \textbf{Periodic Simple-Averaging SGD (PSASGD)}~\cite{mcdonald2010distributed} reduces communication frequency by allowing multiple local updates before synchronization across workers.
    \item \textbf{Elastic Averaging SGD (EASGD)}~\cite{zhang2015deep} introduces a central anchor variable to constrain the divergence among local models while preserving flexibility in local updates.
    \item \textbf{Decentralized Parallel SGD (D-PSGD)}~\cite{nedic2018network} enables peer-to-peer model averaging within local neighborhoods, eliminating the reliance on a central server and improving scalability.
\end{itemize}

While numerous algorithmic variants have been proposed, the theoretical understanding—particularly in terms of convergence guarantees—remains limited. For PSASGD, convergence has been established for both strongly convex and non-convex objectives under bounded gradient assumptions~\cite{stich2018local,yu2019parallel,jiang2018linear}. In contrast, EASGD's analysis is restricted to one-step updates and quadratic loss functions~\cite{zhang2015deep}. D-PSGD has been shown to converge for non-convex losses, though typically without allowing for multiple local updates~\cite{lian2017can,jiang2017collaborative,zeng2018nonconvex}. More recently, Cooperative SGD~\cite{wang2021cooperative} introduced a unified convergence analysis framework that encompasses PSASGD, EASGD, and D-PSGD.

Despite these advancements, most theoretical results focus on vanilla SGD without momentum. However, momentum-based variants are widely adopted in practice due to their superior convergence speed and generalization performance~\cite{Krizhevsky2012ImageNet,yan2018unified,sutskever2013importance}. Among existing analyses, the only known result for momentum-based distributed SGD establishes time-averaged convergence for non-convex objectives~\cite{yu2019linear}. Yet, last-iterate convergence—more relevant to practical training—remains largely unexplored. This reveals a significant gap: although momentum-based distributed methods are predominant in deep learning applications, their theoretical foundations, especially regarding last-iterate convergence, are still underdeveloped.

\textbf{Our Contributions.} 
This paper investigates the momentum-based variants of three widely used distributed optimization methods—PSASGD, EASGD, and D-PSGD—under the classical Robbins–Monro step-size schedule. The main contributions are summarized as follows.

\begin{itemize}
    \item In non-convex settings, We provide the first theoretical analysis that establishes both \textbf{almost sure} and \textbf{mean-square} last-iterate convergence of the gradient norm for such momentum-based algorithms. 
    \item We further derive convergence rate estimates.
    \item We show that momentum accelerates early-stage convergence toward a neighborhood of stationary points. Experimental results support our theoretical findings.
\end{itemize}

To the best of our knowledge, this is the first work to establish last-iterate convergence results for distributed momentum-based PSASGD, EASGD, and D-PSGD algorithms.

\section{Main Results} \label{dpsgdertr}

\subsection{Definitions of Convergence}

For the problem \eqref{distributed_sgd}, suppose the gradient of the loss function \( g_i(x) \) exists and is denoted by \( \nabla g(x) \). Then, for an iterative sequence \( \{x_n\} \), we define the following types of convergence:
\begin{itemize}
   \item \textbf{\(\epsilon\)-TAMS Convergence (\(\epsilon\)-Neighborhood Time-Averaged Mean-Square):} For any scalar \( \epsilon > 0 \), there exists an \( n \) such that
    \[
    \frac{1}{n} \sum_{k=1}^{n} \mathbb{E}\big[\|\nabla g(x_k)\|^2\big] < \epsilon.
    \]

    \item \textbf{TAMS Convergence (Time-Averaged Mean-Square):}
    \begin{multline}
        \frac{1}{n} \sum_{k=1}^{n} \mathbb{E}\left[\|\nabla g(x_k)\|^2\right] = O(f(n)), \\
        \text{where } f(n) \to 0 \text{ as } n \to \infty.
    \end{multline}
    
    \item \textbf{LIMS Convergence (Last-Iterate Mean-Square):}
    \begin{multline}
        \mathbb{E}\left[\|\nabla g(x_n)\|^2\right] = O(f(n)), \\
        \text{where } f(n) \to 0 \text{ as } n \to \infty.
    \end{multline}

    \item \textbf{LIAS Convergence (Last-Iterate Almost-Sure):}
    \begin{multline}
        \|\nabla g(x_n)\| = O(f(n)) \quad \text{almost surely}, \\
        \text{where } f(n) \to 0 \text{ as } n \to \infty.
    \end{multline}
\end{itemize}

We note that LIMS convergence implies TAMS convergence, but the converse does not hold. Similarly, TAMS convergence implies \(\epsilon\)-TAMS convergence.

\subsection{General Momentum-Based Iteration}

We begin by reviewing two existing distributed SGD algorithms.

\paragraph{D-PSGD.} The decentralized SGD algorithm (D-PSGD) has been studied in various works \cite{jiang2021consensus,2017Deep,2017Can,wang2021cooperative}. In this approach, each worker node performs local updates and synchronizes with neighbors every \( k \) iterations:
\[
\begin{aligned}
    &\text{if } n \bmod k = 0: \\
    &\quad x_{n+1}^{(i)} = \sum_{j=1}^{m} w_{ji} \left(x_n^{(j)} - \epsilon_n \nabla g_j(x_n^{(j)}, \xi_n^{(j)}) \right), \\
    &\text{else:} \\
    &\quad x_{n+1}^{(i)} = x_n^{(i)} - \epsilon_n \nabla g_i(x_n^{(i)}, \xi_n^{(i)}),
\end{aligned}
\]
where the entry \( w_{ji} \) of the mixing matrix \( W \), located at the \((j,i)\)-th position, quantifies the influence of node \( j \) on the model averaged at node \( i \).
 \textbf{PSASGD} is a special case of D-PSGD with \( w_{ji} = \frac{1}{m} \).

\paragraph{EASGD.} In contrast, EASGD uses an elastic force mechanism to coordinate local models, inspired by quadratic penalty methods \cite{zhang2015deep}:
\begin{equation} \label{EASGD2}
\begin{aligned}
    &\text{if } n \bmod k = 0: \\
    &\quad x_{n+1}^{(i)} = (1 - \beta)\left(x_n^{(i)} - \epsilon_n \nabla g_i(x_n^{(i)}, \xi_n^{(i)})\right) + \beta z_n, \\
    &\quad z_{n+1} = (1 - m\beta)z_n + m\beta \bar{x}_n, \\
    &\text{else:} \\
    &\quad x_{n+1}^{(i)} = x_n^{(i)} - \epsilon_n \nabla g_i(x_n^{(i)}, \xi_n^{(i)}), \\
    &\quad z_{n+1} = z_n,
\end{aligned}
\end{equation}
where \( \bar{x}_n = \frac{1}{m} \sum_{i=1}^m x_n^{(i)} \), and \( \beta > 0 \) is the consensus control parameter.

A unified update rule for EASGD and D-PSGD is given by:
\begin{equation} \label{combine}
    X_{n+1} = W_n \left(X_n - \epsilon_n G(X_n, \xi_n) \right),
\end{equation}
with:

- For D-PSGD:
\[
\begin{aligned}
    X_n &= (x_n^{(1)}, x_n^{(2)}, \dots, x_n^{(m)})^\top, \\
    G(X_n, \xi_n) &= \left(\nabla g_1(x_n^{(1)}, \xi_n^{(1)}), \dots, \nabla g_m(x_n^{(m)}, \xi_n^{(m)})\right)^\top, \\
    W_n &= \begin{cases}
        (w_{ij})_{m \times m}, & \text{if } n \bmod k = 0, \\
        I_m, & \text{otherwise}.
    \end{cases}
\end{aligned}
\]
The PSASGD algorithm is a special case of the D-PSGD algorithm when \( W = \frac{\mathbf{1} \mathbf{1}^\top}{\mathbf{1}^\top \mathbf{1}} \), where \( \mathbf{1} = (1, 1, \dots, 1)^\top \) is a column vector of size \( m \). It will therefore not be discussed separately in the remainder of the paper.

- For EASGD:
\begin{equation}
\begin{aligned}
    X_n &= (x_n^{(1)}, \dots, x_n^{(m)}, z_n)^\top, \\
    G(X_n, \xi_n) &= \left(\nabla g_1(x_n^{(1)}, \xi_n^{(1)}), \dots, \nabla g_m(x_n^{(m)}, \xi_n^{(m)}), \mathbf{0}\right)^\top, \\
    W_n &= 
    \begin{cases}
        \begin{pmatrix}
            (1{-}\beta)I & \beta \mathbf{1} \\
            \beta \mathbf{1}^\top & 1 {-} m\beta
        \end{pmatrix}, & \text{if } n \bmod k = 0, \\[5pt]
        I, & \text{otherwise}.
    \end{cases}
\end{aligned}
\end{equation}
The term \(\mathbf{0}\) denotes a \(v \times 1\) zero vector.

\paragraph{Momentum Extension.} To accelerate convergence, we incorporate momentum into \eqref{combine}, yielding:
\begin{equation} \label{123combinemomentum}
\begin{aligned}
    v_n &= \alpha v_{n-1} + \epsilon_n G(X_n, \xi_n), \\
    X_{n+1} &= W_n \left(X_n - v_n \right),
\end{aligned}
\end{equation}
where \( \alpha \in [0,1) \) is the momentum coefficient, \( \epsilon_n \) is the step size, and \( v_0 := 0 \). 
This general form includes the momentum-based EASGD \cite{zhang2015deep}, D-PSGD \cite{yu2019linear, 2021DecentLaM,2020SQuARM,2020Periodic,2021Decentralized} as a special case.

Let \( X = (x^{(1)}, \dots, x^{(m)}) \). Then:\\
- For D-PSGD:
\[
\begin{aligned}
    G(X, \xi_n) &= \left(\nabla g_1(x^{(1)}, \xi_n^{(1)}), \dots, \nabla g_m(x^{(m)}, \xi_n^{(m)})\right)^\top, \\
    G(X) &= \left(\nabla g_1(x^{(1)}), \dots, \nabla g_m(x^{(m)})\right)^\top.
\end{aligned}
\]
- For EASGD:
\[
\begin{aligned}
    G(X, \xi_n) &= \left(\nabla g_1(x^{(1)}, \xi_n^{(1)}), \dots, \nabla g_m(x^{(m)}, \xi_n^{(m)}), \mathbf{0}\right)^\top, \\
    G(X) &= \left(\nabla g_1(x^{(1)}), \dots, \nabla g_m(x^{(m)}), \mathbf{0} \right)^\top.
\end{aligned}
\]

In the following sections, we study the convergence properties of the general momentum-based iteration \eqref{123combinemomentum}.

\subsection{Last-iterate Convergence}
To proceed, the following assumptions are needed. 
\begin{assum}\label{ass_g1} $g(x):=\frac{1}{m}\sum_{i=1}^{m}g_{i}(x)$ is a non-negative and  continuously differentiable. In addition, the following conditions hold:
	\begin{enumerate}
		\item $G(X,\xi_{n}))$ is an unbiased estimate of $G(X)$, i.e., $\Expect_{\xi_{n}}G(X,\xi_{n})=G(X)$;
		\item  The mixing matrix $W_{n}\in \mathbb{R}^{m\times m}$ is a symmetric doubly stochastic matrix with only one eigenvalue equal to one and the absolute values of the rest eigenvalues are less than one. 
		\item (Assumption 1 in \cite{yu2019linear}) There are two constants $L>0,\ M>0$, such that $\forall, X, Y\in\mathbb{R}^{m\times N}$, $\|G(X)-G(Y)\|\le L\|X-Y\|$ and $\|G(X)\|\le M.$ 
		
		\item For any  $ i=1,2,...,m$ and $\forall X\in\mathbb{R}^{m\times N}$, it holds that 
  \begin{equation}\nonumber\begin{aligned}
				&\sum_{i=1}^{m}\Expect_{\xi_{n}}\big\|\nabla g_{i}(x,\xi)-\nabla g_{i}(x)\big\|^{2}\le\sigma_{0}^{2}\ .\end{aligned}\end{equation}

  In addition, $\forall \ x\in\mathbb{R}^{N}$, it holds that
		\begin{equation}\nonumber\begin{aligned}
				&\frac{1}{m}\sum_{i=1}^{m}\|\nabla g_{i}(x)-\nabla g(x)\|^{2}\le \sigma_{1}^{2}\ .
		\end{aligned}\end{equation}
	\end{enumerate}
\end{assum}
The conditions in assumption \ref{ass_g1}  are common in the study of   distributed SGD or mSGD. We can find these conditions in the literature \cite{yu2019linear,wang2021cooperative,2019Parallel,jin2022convergence,Nguyen2018SGD}. In some works, the non-negative loss function condition may be replaced by a finite low bound condition, i.e., $g(x)>\hat{l}_{0}>-\infty$. These two conditions are essentially equivalent, since  one can construct a new loss function $\overline{g}=g-\hat{l}_{0}$ for the finite low bound condition, such that the new loss function is non-negative. Note that item 4 in Assumption 1 quantifies the variance of stochastic gradients at local worker, and $\sigma_{1}^2$ quantifies the deviations between the local objective function of each workers.



\begin{assum}\label{opjknhg} 
	The momentum coefficient $\alpha\in [0,1)$ and 
	the sequence of  learning rate $\epsilon_{n}$ satisfies the Robbins-Monro condition, i.e., it is positive,  monotonically decreasing to zero, such that
	$	\sum_{n=1}^{+\infty}\epsilon_{n}=+\infty$ and 
	$\sum_{n=1}^{+\infty}\epsilon_{n}^{2}<+\infty.$
\end{assum}

As stated in Assumption \ref{opjknhg}, we focus on the optimization process with a decaying learning rate, as fixed step sizes may cause the model to oscillate near the bottom of the loss landscape, preventing stable convergence. In contrast, a decaying schedule helps stabilize training and facilitates convergence \cite{smith2017don,2011Bayesian,2015Convergence,2019Understanding,2016Deep,jin2022convergence}.

Under the above assumptions, we establish the convergence of momentum-based distributed SGD iterates as described in Equation \eqref{123combinemomentum}, which is formalized in the following theorem.
\begin{thm}\label{lem_mied}
	Suppose $\{X_{n}\}$ is a sequence  generated by equation \eqref{123combinemomentum}. Under Assumptions \ref{ass_g1}--\ref{opjknhg}, it holds that $\|\nabla  g(\overline{x}_{n})\|\rightarrow 0\ a.s.$ and $\Expect\|\nabla  g(\overline{x}_{n})\|^{2}\rightarrow 0,$ where $\overline{x}_{n}$ is defined as the average value of every worker node, i.e., $\overline{x}_{n}={1}/{m}\sum_{i=1}^{m}x_{n}^{(i)}\ .$
\end{thm}

Theorem~\ref{lem_mied} builds upon the centralized momentum-based SGD framework proposed in~\cite{jin2022convergence}, introducing several key innovations to extend its applicability to distributed settings. It establishes the last-iterate convergence of the iterates in Equation~\eqref{123combinemomentum}, and this result also implies convergence in the time-averaged sense, as considered in~\cite{yu2019linear,2021DecentLaM,2020SQuARM,2020Periodic,2021Decentralized}; that is,
\[
\frac{1}{T} \sum_{i=1}^{T} \mathbb{E}\left[\|\nabla g(\overline{x}_{n})\|^{2}\right] \rightarrow 0.
\]

\subsection{Last-iterate convergence rate}
To quantitatively estimate the convergence rate of the last iterate, an additional Assumptions \ref{ass_steps123} and \ref{ass_steps1234} are required in this proof. In existing literature, the strong convexity assumption is commonly adopted. For example, \cite{2021DecentLaM} assumes that each local loss function \( g_i \) is strongly convex when analyzing the last-iterate convergence rate of deterministic distributed momentum-based gradient descent. Similarly, \cite{Nguyen2018SGD} requires the strong convexity of the overall loss function \( g \) when studying the last-iterate convergence rate of SGD. 
\begin{assum}\label{ass_steps123} The loss function \(g(\theta)\) is a convex function and has a unique optimal point \(\theta^*\).
\end{assum}
\begin{assum}\label{ass_steps1234}
During the algorithm iteration process, stability is maintained, i.e., for any \(n > 0\), there exists a constant \(G < +\infty\) such that \(\|u^{\top}X_{n}\| < G\) almost surely.
\end{assum}

Under these additional assumptions, we establish the last-iterate convergence rate as follows:

\begin{thm}\label{coro_41}
	Suppose $\{X_{n}\}$ is a sequence  generated by equation \eqref{123combinemomentum}. Under Assumptions \ref{ass_g1}-- \ref{ass_steps1234} with $\epsilon_{n}=\frac{\sqrt{m}}{\sqrt{n}}.$ Then for any $T>0,$ there is
	\[\Expect(g(u^{\top}X_{T})-g(\theta^*))=\mathcal{O}\Big(\sqrt{m}\frac{\ln T}{\sqrt{T}}\Big)+\mathcal{O}\Big(\frac{1}{\sqrt{m}}\frac{\ln T}{\sqrt{T}}\Big).\] 
\end{thm}

It can be observed that the incorporating of momentum does not significantly improve the algorithm's asymptotic convergence rate. This observation appears inconsistent with empirical results, which show that momentum can substantially accelerate convergence in practice. The root of this discrepancy lies in the fact that the convergence rate discussed here characterizes asymptotic behavior as the number of epochs approaches infinity, whereas momentum primarily accelerates convergence in the early stages of training. To formalize this phenomenon, we present the following theorem.

\begin{thm}\label{dacnie}
	Suppose $\{X_{n}\}$ is a sequence generated by equation \eqref{123combinemomentum}. Under Assumption \ref{ass_g1}, given any non-increasing positive learning rate $\epsilon_{n}\ge \epsilon_{n+1}$ and bounded loss function,  for any worker node $i\ \ (i=1,2,...,m)$,  then for any $a_{0}>0$,  any $ V_{0}\in \mathbb{R}^{mN}$ and any $ \  \|\nabla  g(\overline{x}_{1})\|^{2}>a_{0}$, there exists $s>0$, such that
	\begin{equation}\nonumber\begin{aligned}
			&P({\tau}^{(a_{0})}\ge n)=O\Big(e^{-\frac{s}{(1-\alpha)^{2}}\sum_{i=1}^{n}\epsilon_{i}}\Big),
	\end{aligned}\end{equation}  where  ${\tau}^{(a_{0})}=\min_{n>0}\{\|\nabla  g_{i}(x_{n})\|^{2}<a_{0}\}$.
\end{thm}

Theorem \ref{dacnie} shows that  a larger momentum term coefficient $\alpha$  can speed up the convergence in an early stage. 
More specifically, for any small constant \( a_0 \) and iteration number \( n \), the larger the momentum coefficient, the smaller the probability that the number of iterations required for the gradient norm to first fall below \( a_0 \) exceeds \( n \).

\section{Experiments Results}
In this section, we consider a classification task where neural networks are trained using a distributed mSGD algorithm, to demonstrate the correctness of our theoretical findings.

\begin{figure*}[!t]
	\centering
	\subfloat[]{
		\includegraphics[width=0.24\textwidth,height=0.1\textheight]{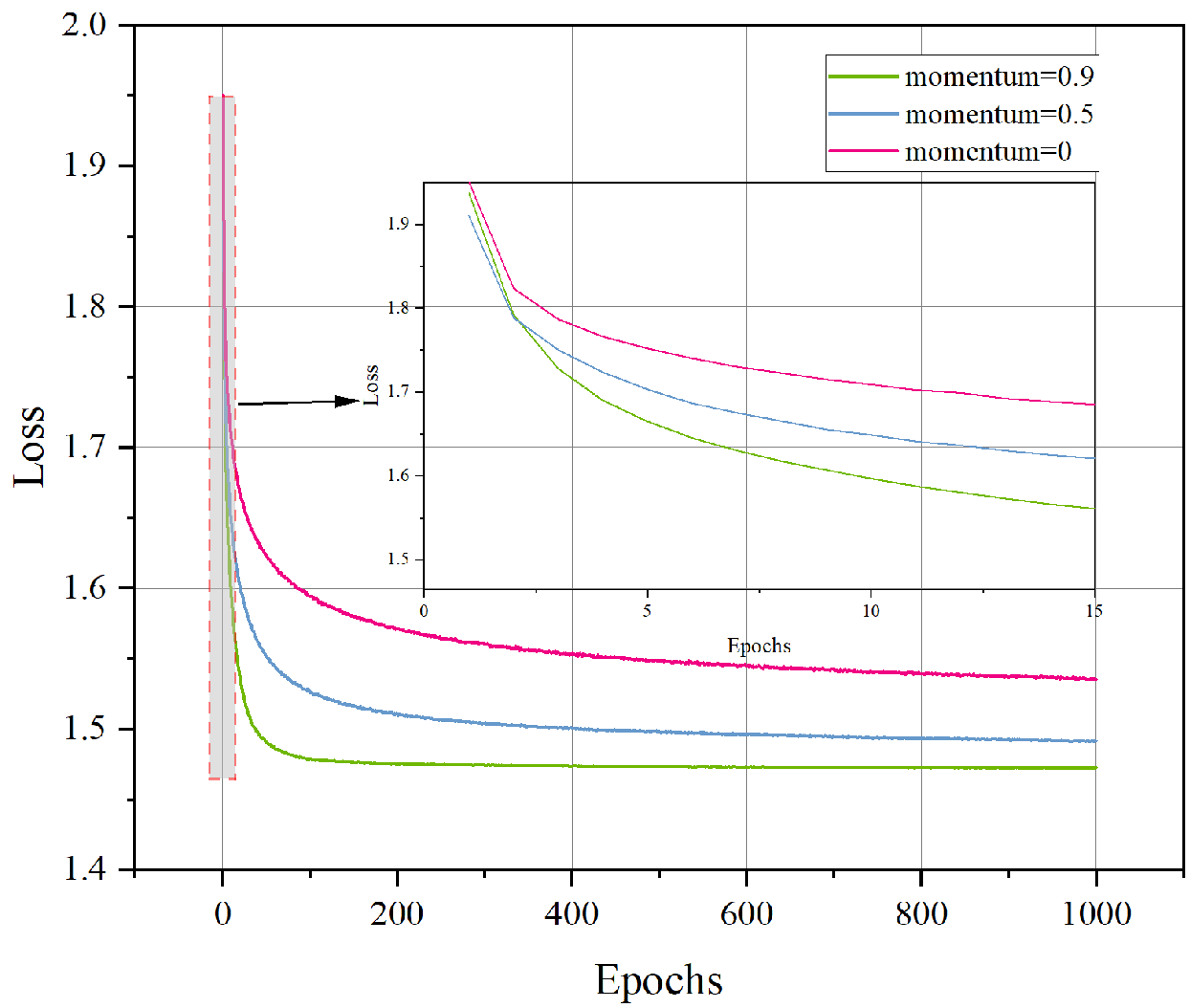}\label{noise_free_10_1}}
  	\hfil
	\subfloat[]{
		\includegraphics[width=0.24\textwidth,height=0.1\textheight]{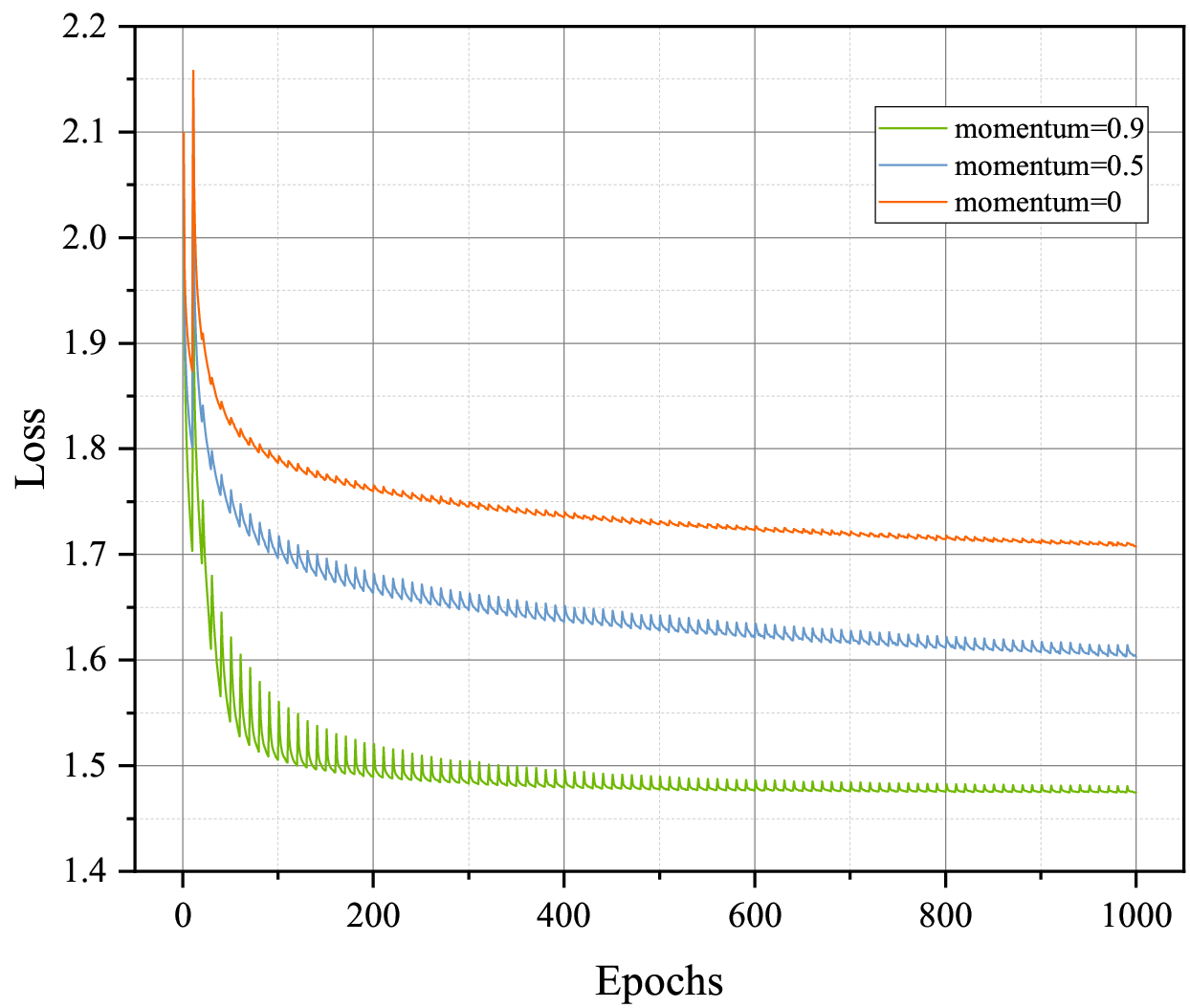}\label{noise_free_10_3}}
  	\hfil
	\subfloat[]{
		\includegraphics[width=0.24\textwidth,height=0.1\textheight]{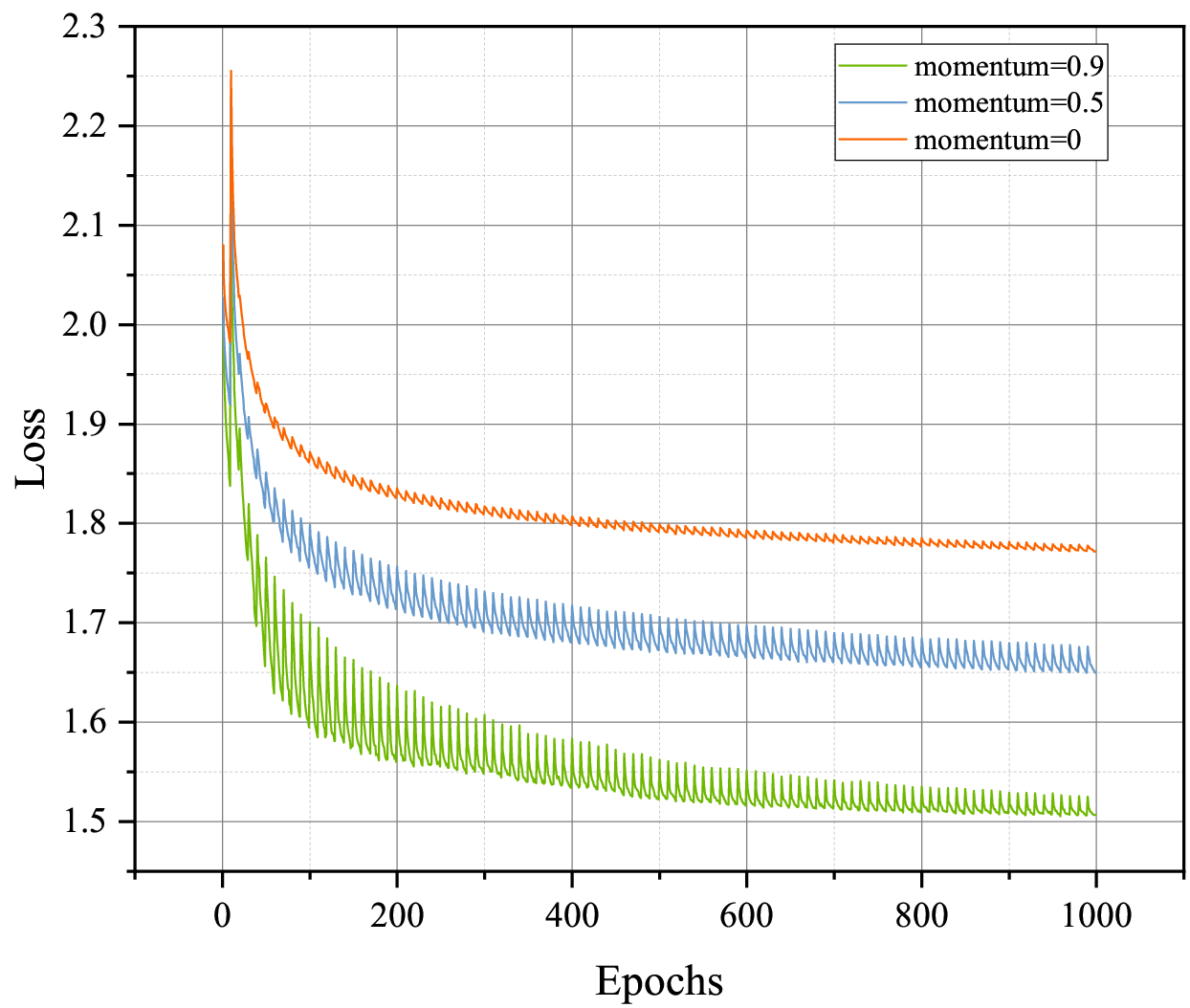}\label{noise_free_10_10}}
  	\hfil
	\subfloat[]{
		\includegraphics[width=0.24\textwidth,height=0.1\textheight]{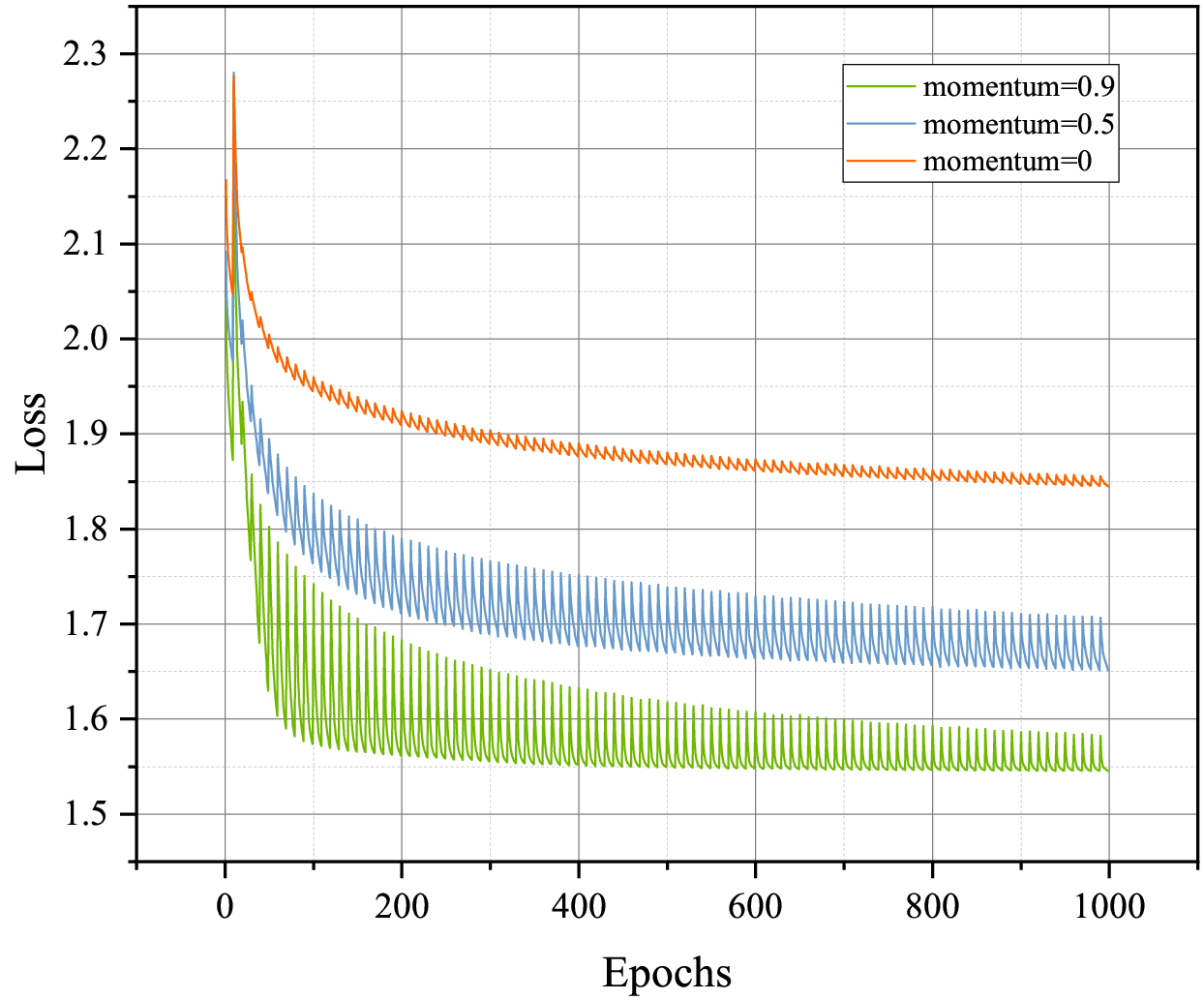}\label{noise_free_10_20}}
	\hfil
	\subfloat[]{
		\includegraphics[width=0.24\textwidth,height=0.1\textheight]{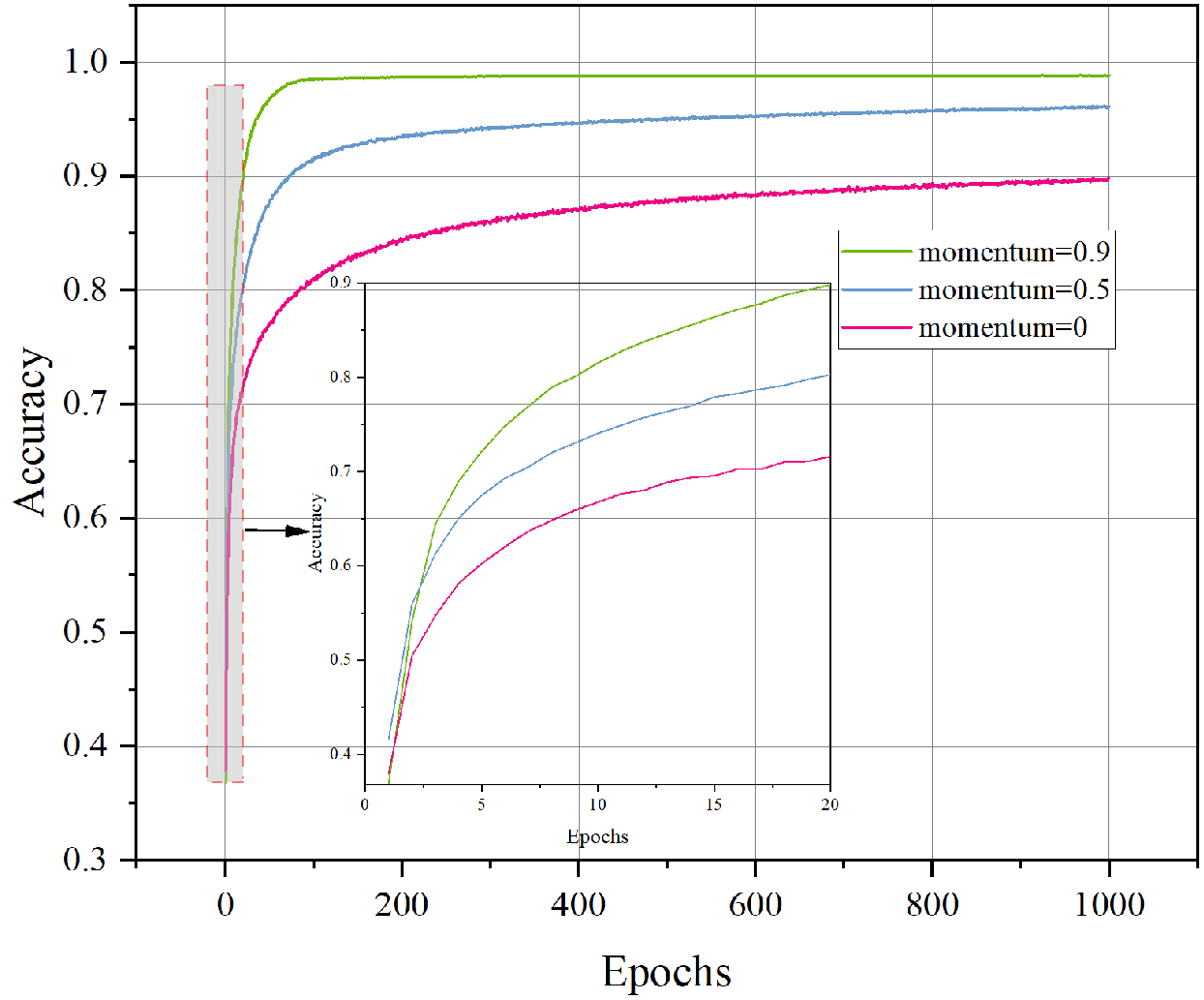}\label{attack_type_10_1}}
	\hfil
	\subfloat[]{
		\includegraphics[width=0.24\textwidth,height=0.1\textheight]{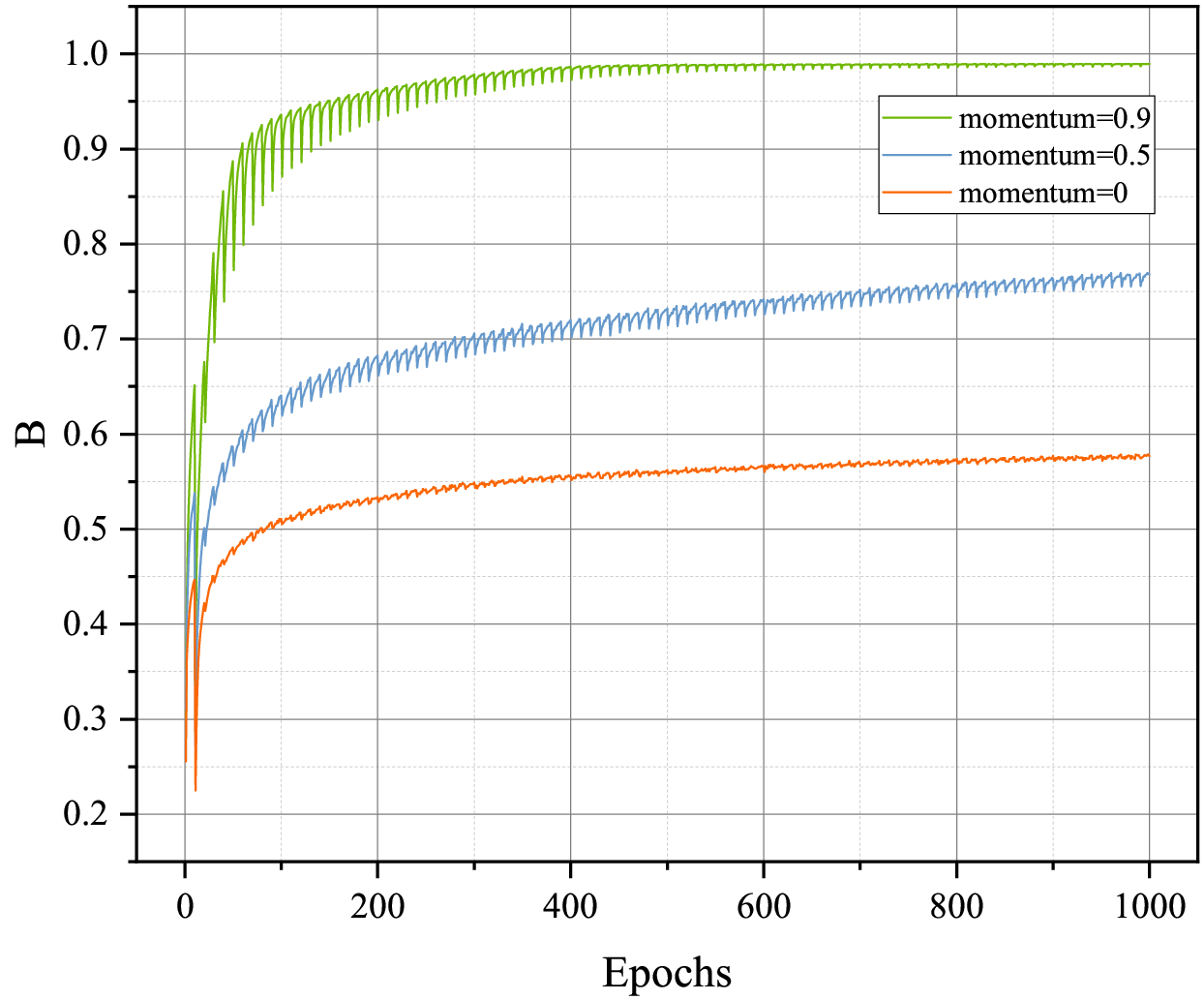}\label{attack_type_10_3}}
	\hfil
	\subfloat[]{
		\includegraphics[width=0.24\textwidth,height=0.1\textheight]{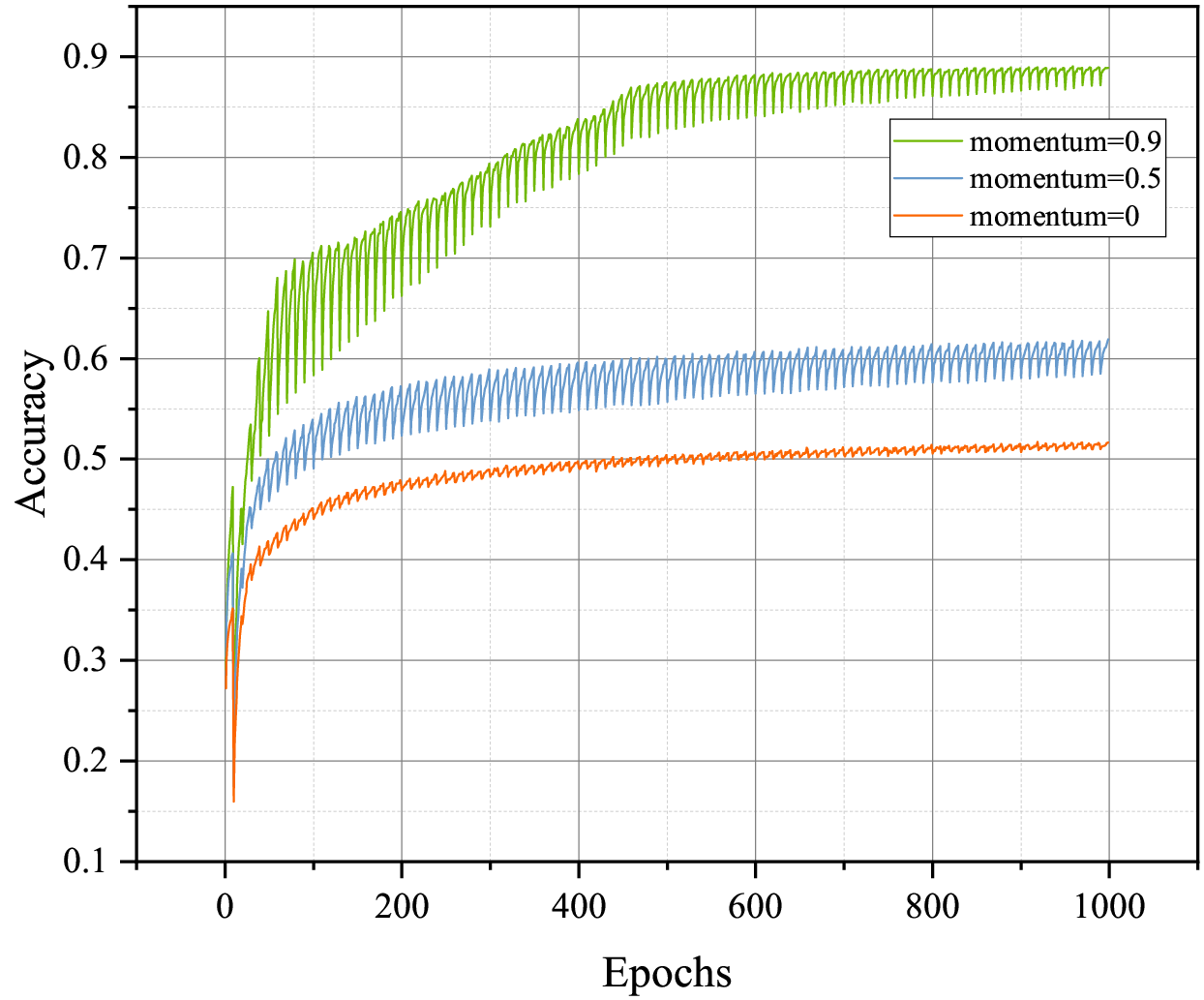}\label{attack_type_10_10}}
	\hfil
	\subfloat[]{
		\includegraphics[width=0.24\textwidth,height=0.1\textheight]{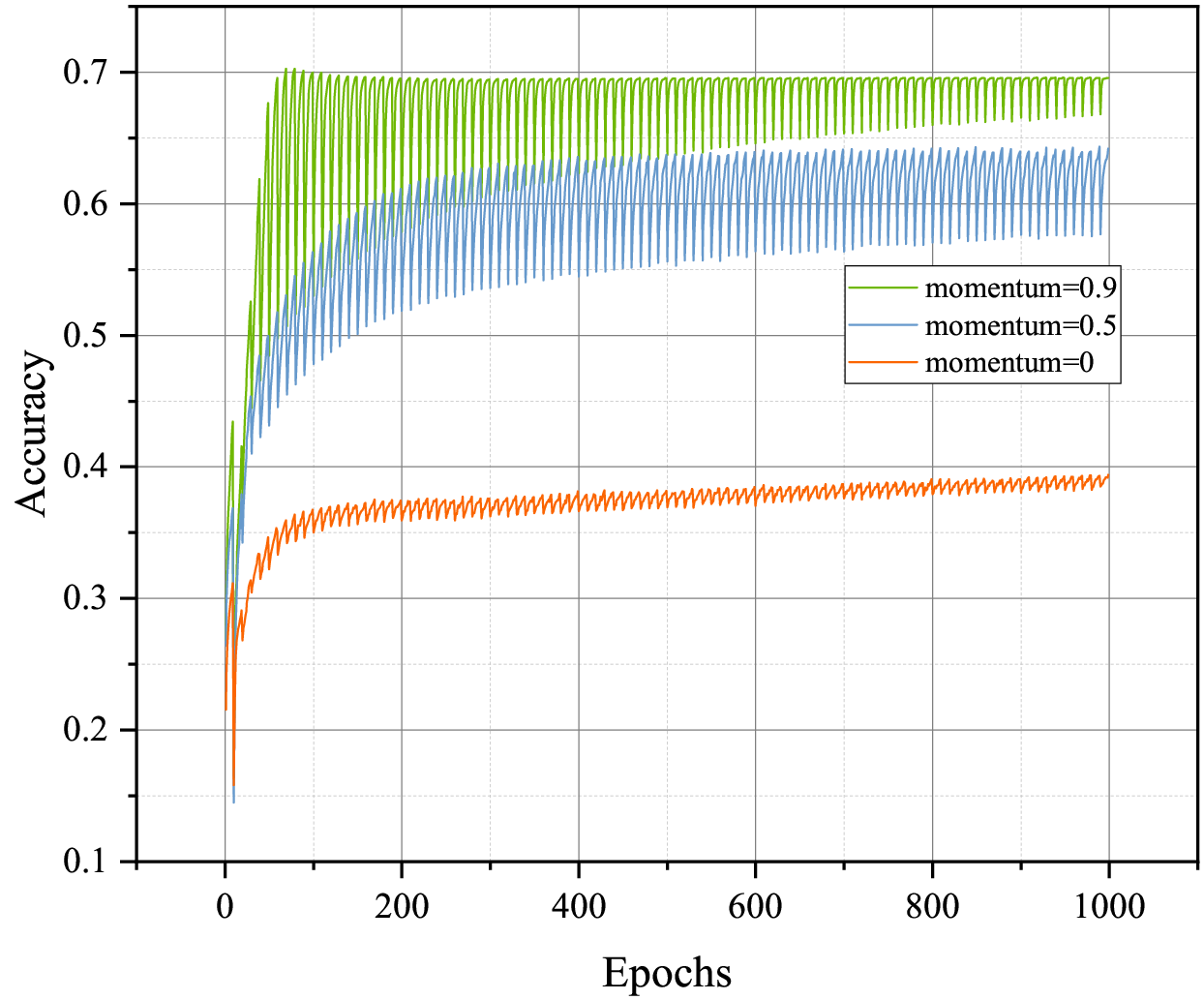}\label{attack_type_10_20}}
    \hfil
    \subfloat[]{
		\includegraphics[width=0.24\textwidth,height=0.1\textheight]{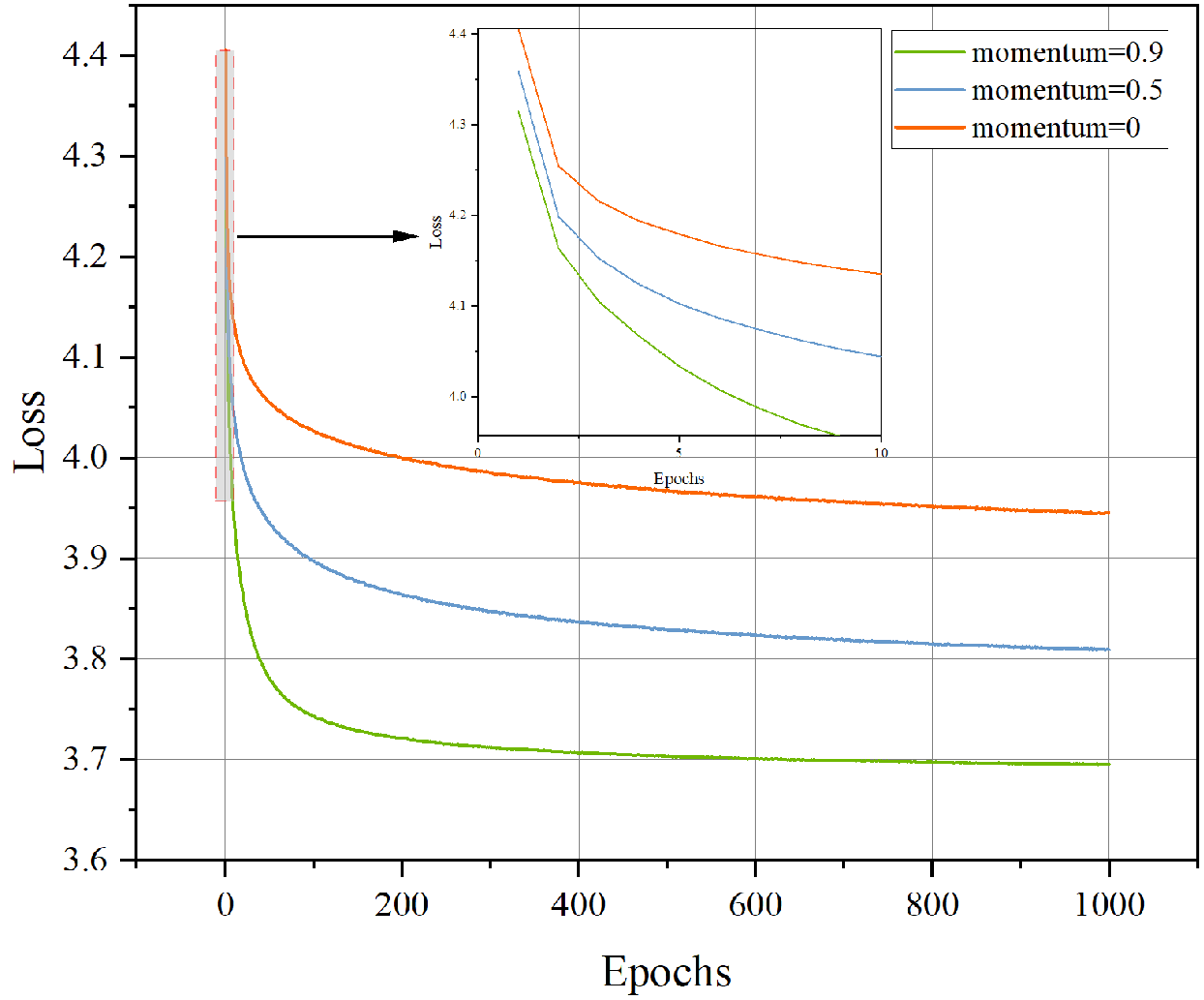}\label{noise_free_100_1}}
        \hfil
	\subfloat[]{
		\includegraphics[width=0.24\textwidth,height=0.1\textheight]{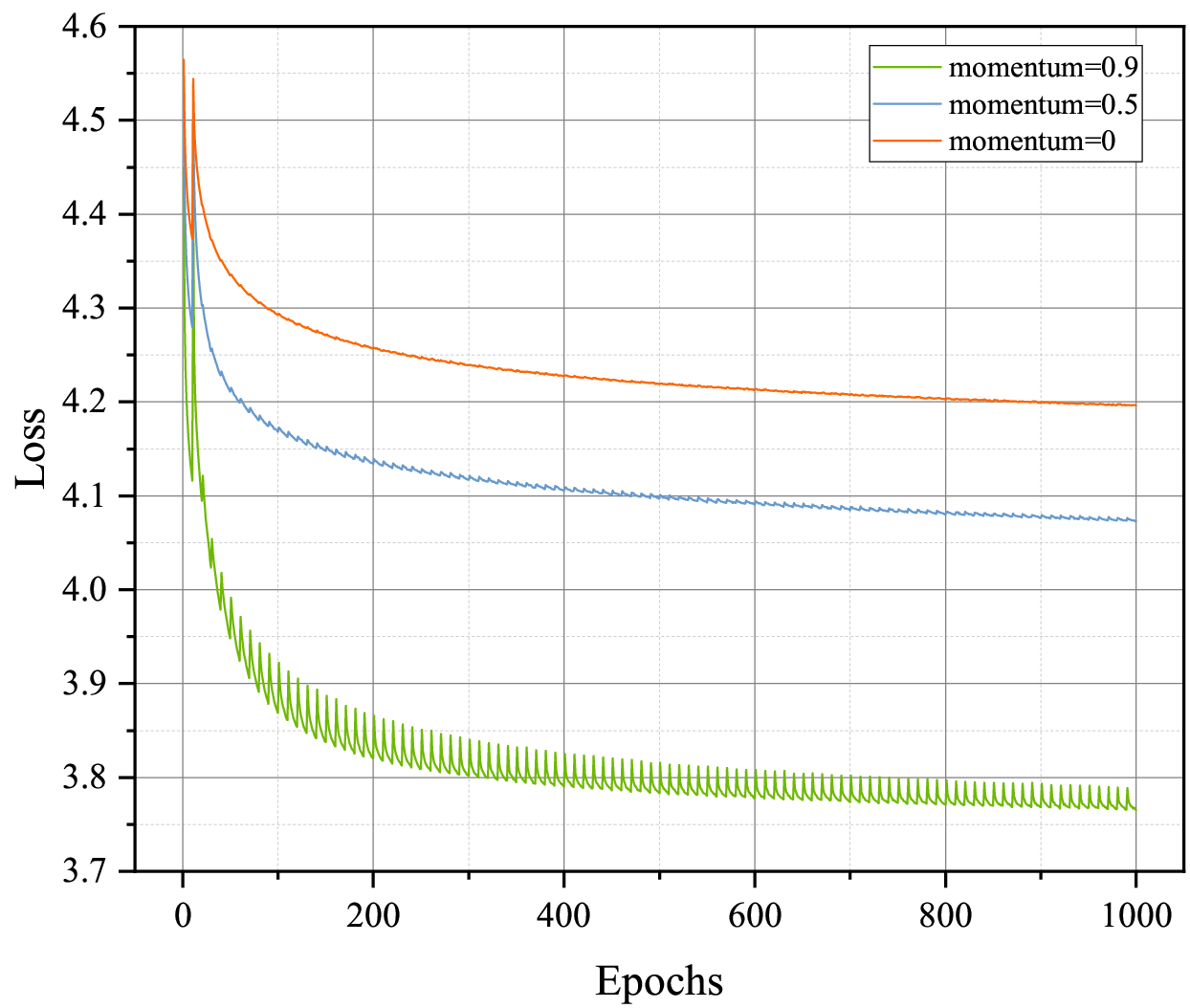}\label{noise_free_100_3}}
  	\hfil
	\subfloat[]{
		\includegraphics[width=0.24\textwidth,height=0.1\textheight]{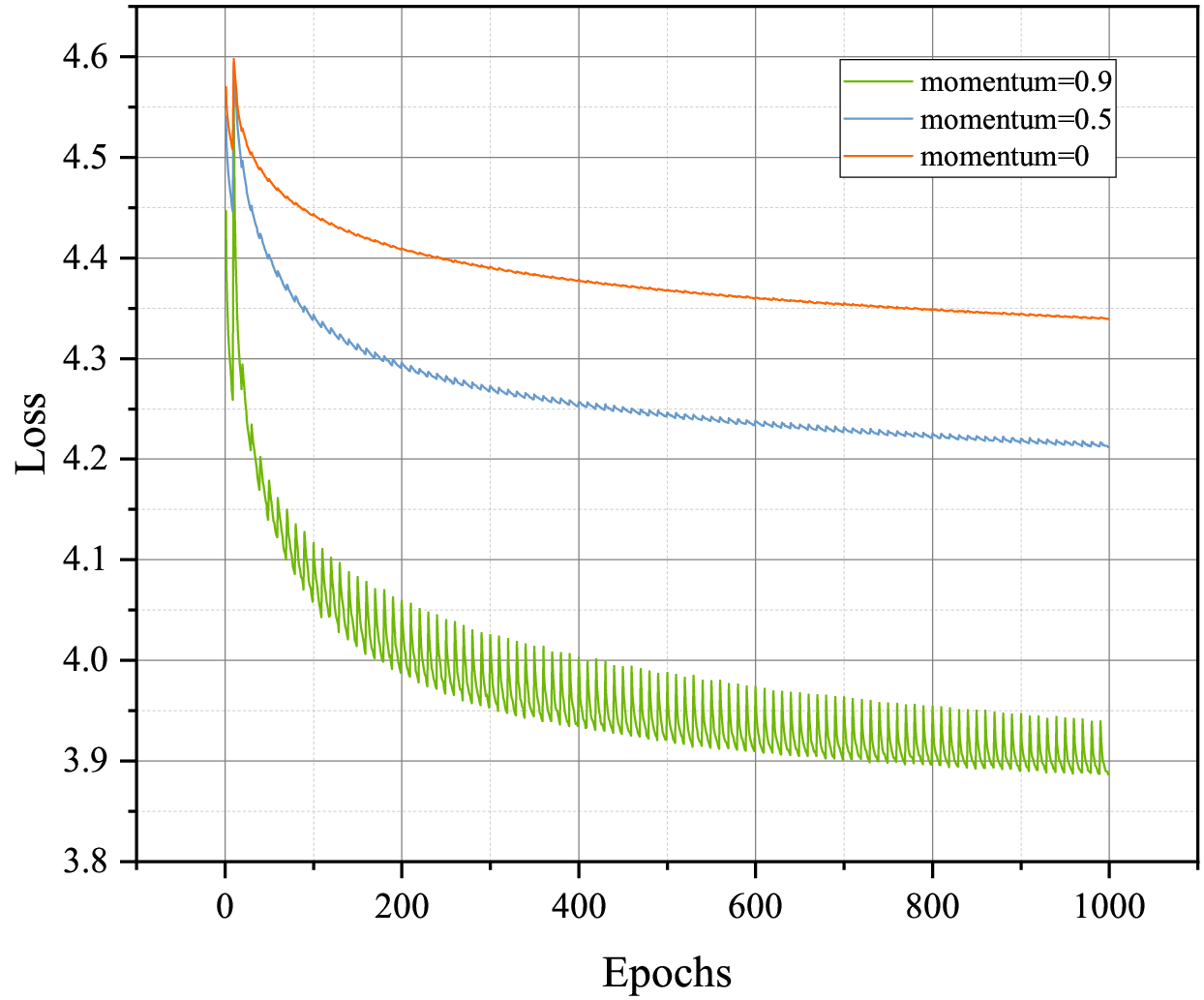}\label{noise_free_100_10}}
  	\hfil
	\subfloat[]{
		\includegraphics[width=0.24\textwidth,height=0.1\textheight]{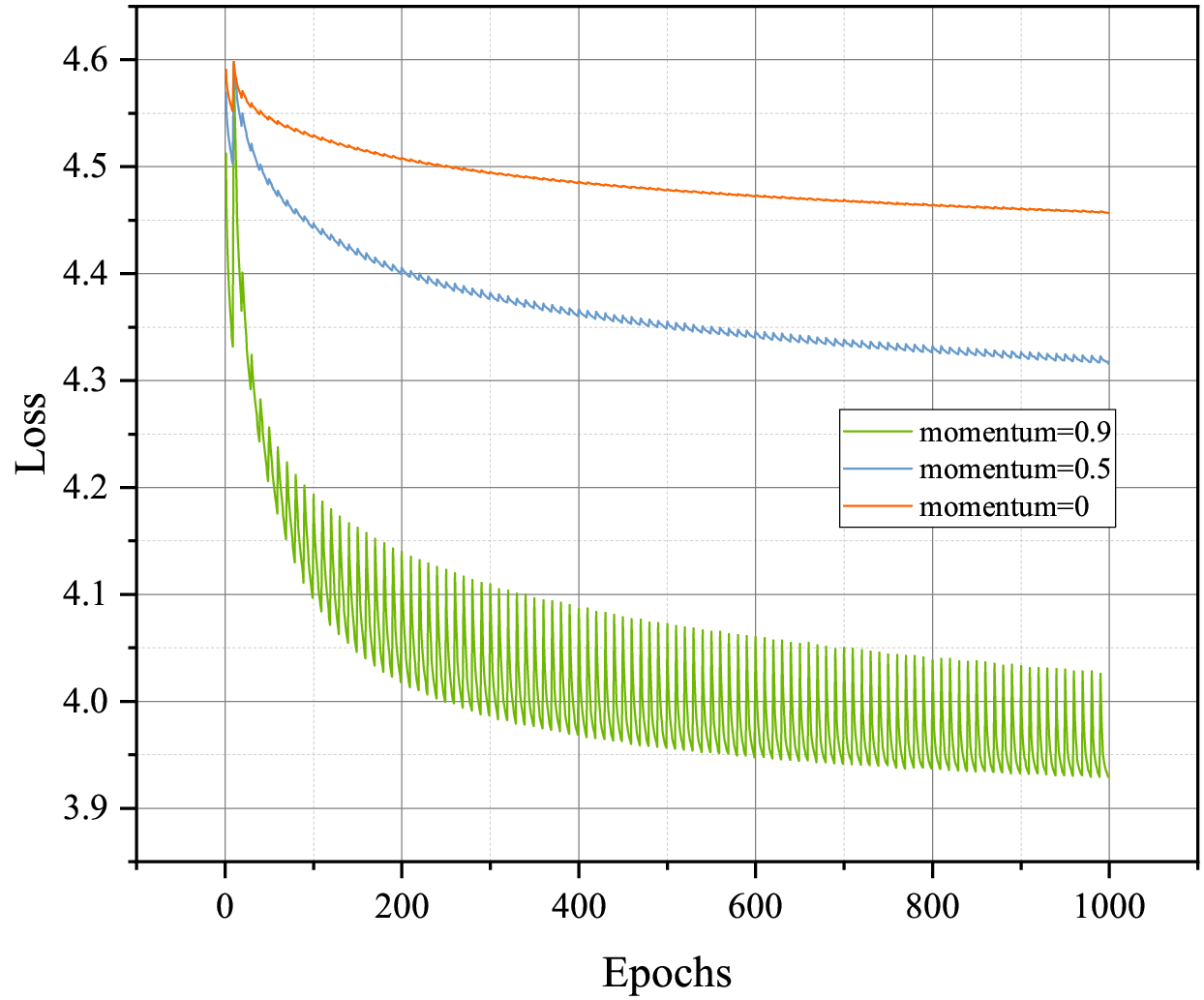}\label{noise_free_100_20}}
	\hfil
	\subfloat[]{
		\includegraphics[width=0.24\textwidth,height=0.1\textheight]{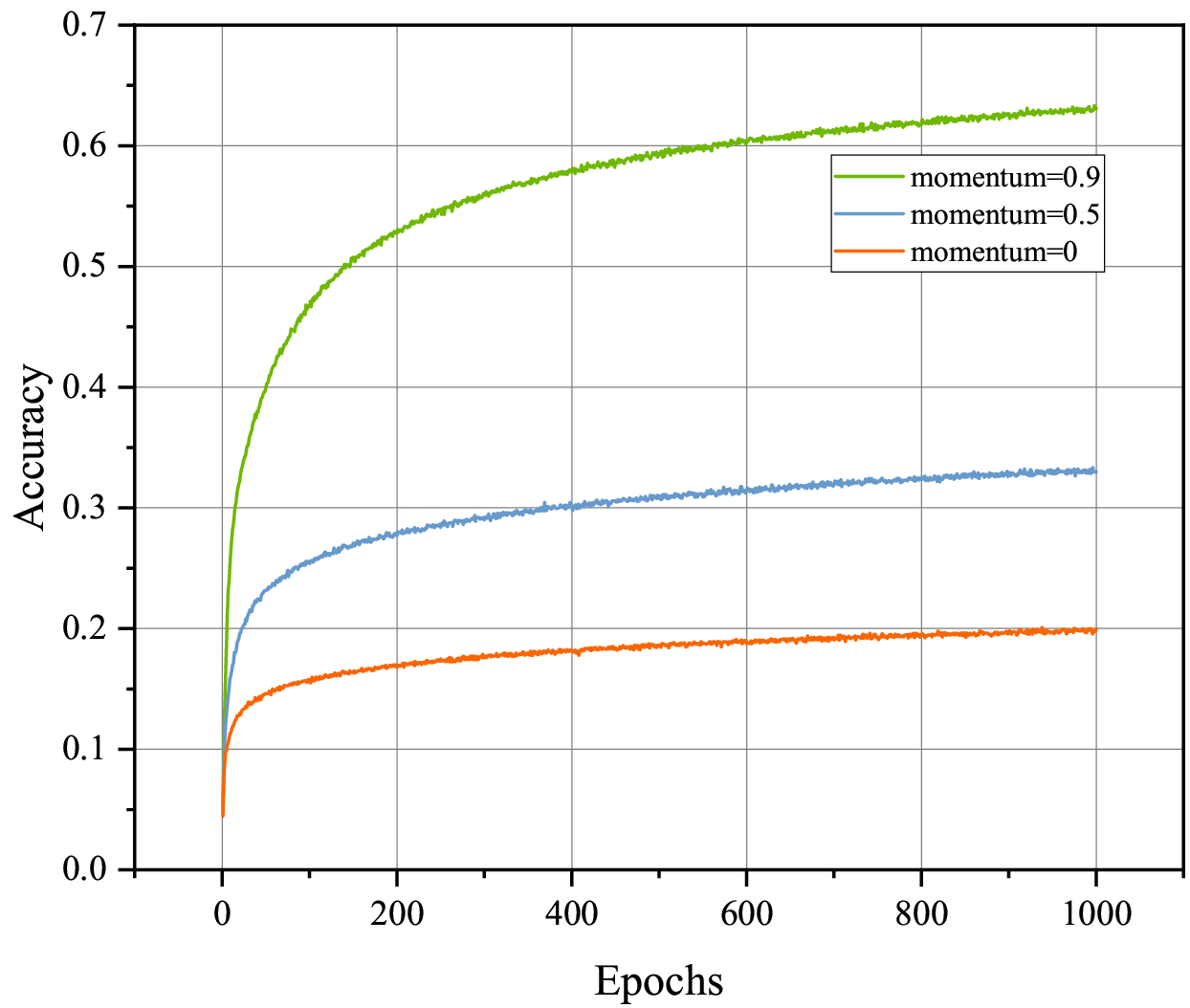}\label{attack_type_100_1}}	
	\hfil
	\subfloat[]{
		\includegraphics[width=0.24\textwidth,height=0.1\textheight]{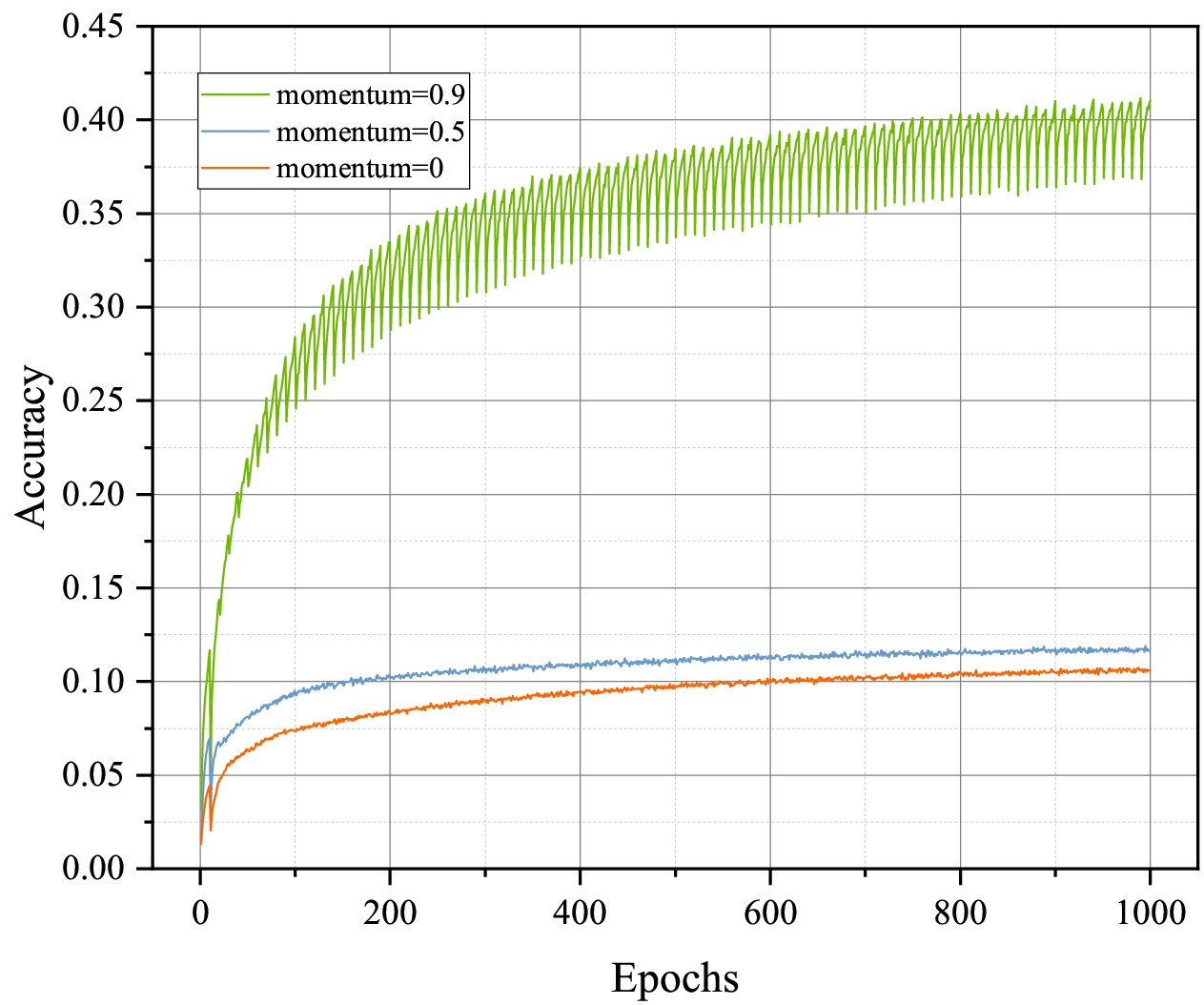}\label{attack_type_100_3}}
	\hfil
	\subfloat[]{
		\includegraphics[width=0.24\textwidth,height=0.1\textheight]{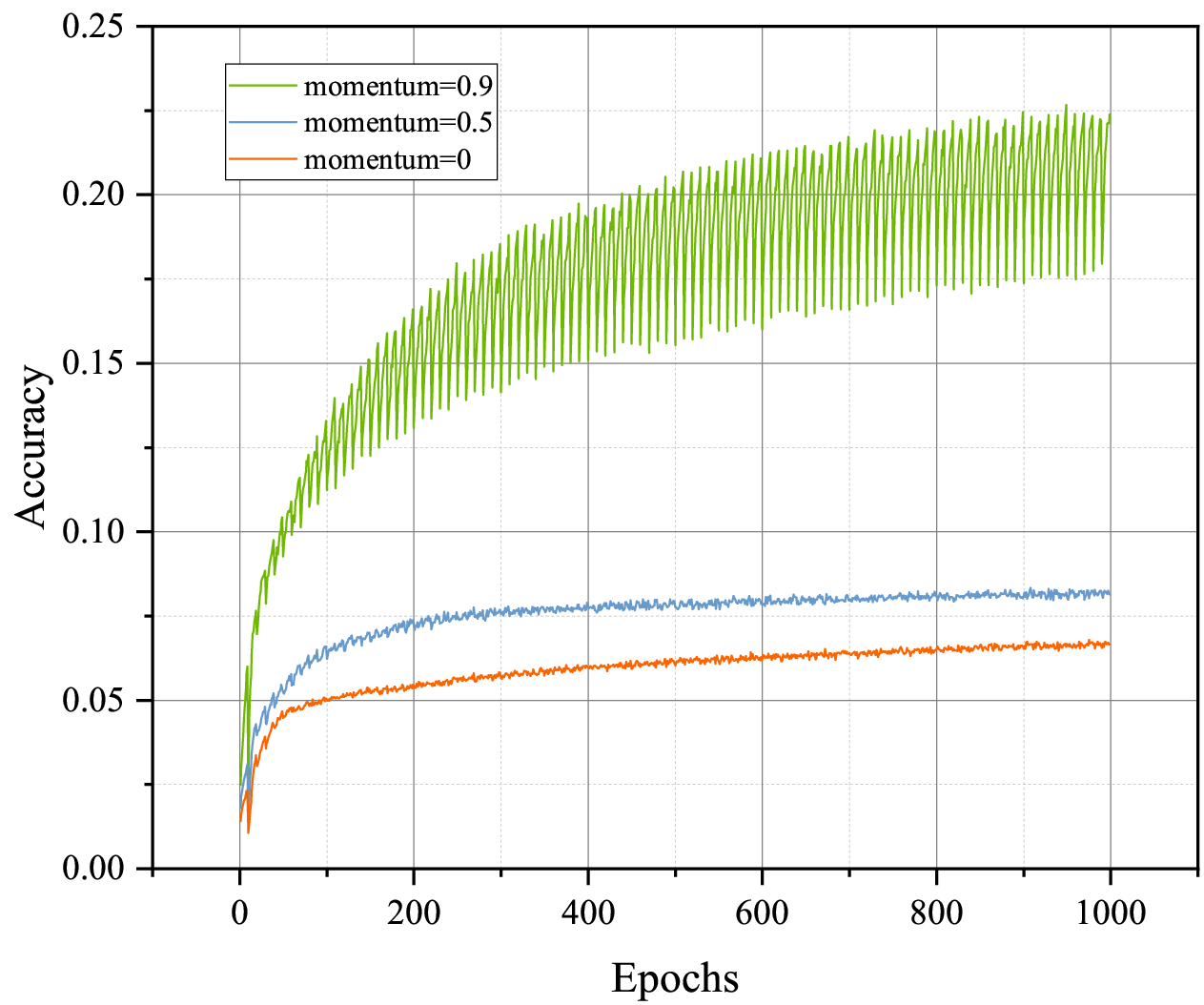}\label{attack_type_100_10}}
	\hfil
	\subfloat[]{
		\includegraphics[width=0.24\textwidth,height=0.1\textheight]{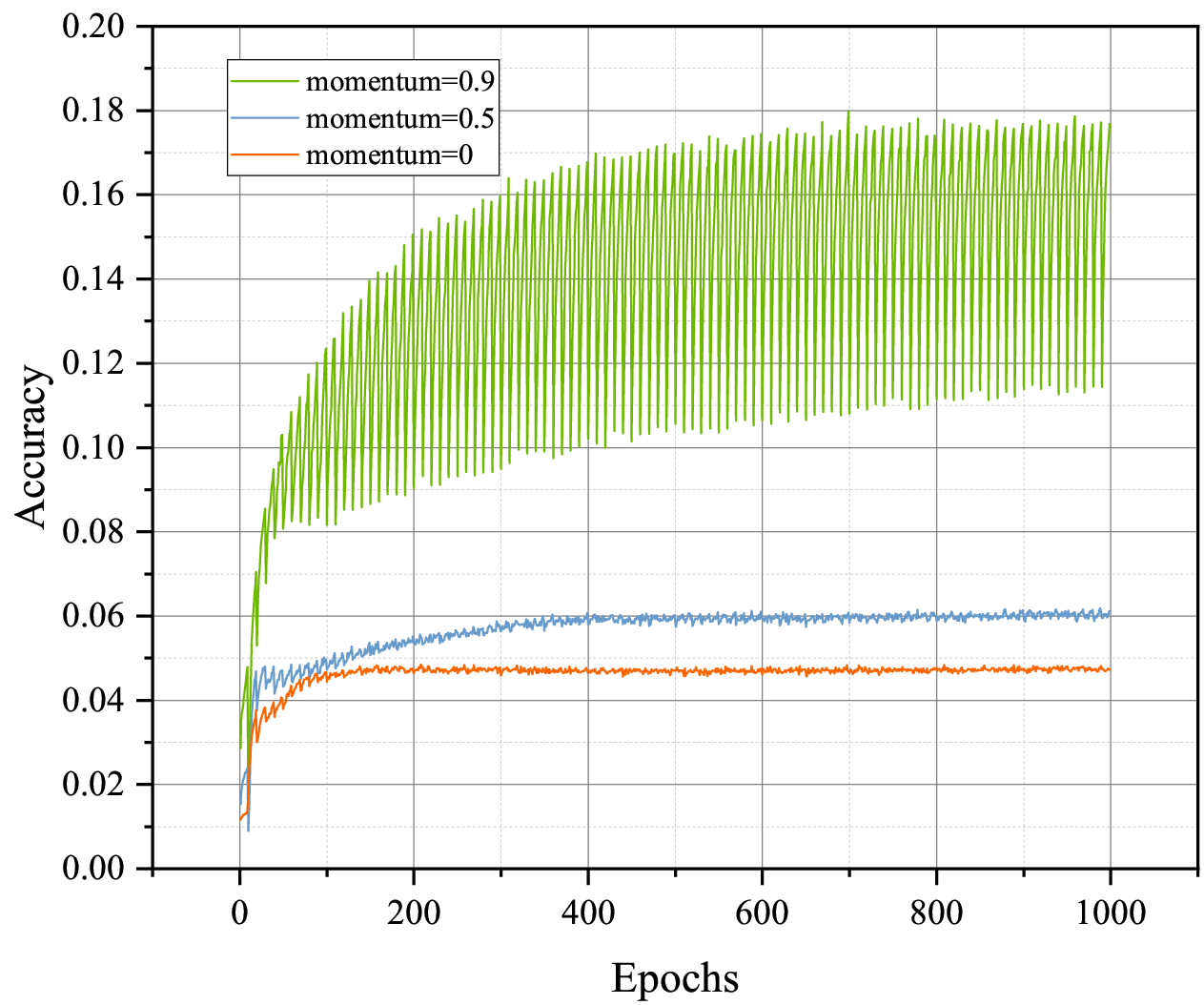}\label{attack_type_100_20}}
	\caption{Training and prediction performance on CIFAR-10 and CIFAR-100 with 1,3,10,20 sub-datasets (workers). (a)-(d): The training loss with 1, 3, 10, and 20 sub-datasets respectively on CIFAR-10. (e)-(h): The accuracy with 1, 3, 10, and 20 sub-datasets respectively on CIFAR-10. (i)-(l): The training loss with 1, 3, 10, and 20 sub-datasets respectively on  CIFAR-100. (m)-(p): The accuracy with 1, 3, 10, and 20 sub-datasets respectively on  CIFAR-100.
	}\label{fig:test1}
\end{figure*}
\textbf{Implementation}.
We employ the ResNet20 network using Keras. We initialize the weights using the Glorot uniform algorithm.
The momentum coefficient takes on the values of 0, 0.5, and 0.9. We train the model using the categorical cross-entropy loss function. The learning rate begins at 0.1 and subsequently decays. We partition the dataset into three, ten, and twenty sub-datasets, with each sub-dataset communicating every 10 epochs with matrices $W$ defined as follows:
$W=\frac{1}{3}\textbf{1}_{3}\textbf{1}_{3}^{\top}$. $W=\frac{1}{10}\textbf{1}_{10}\textbf{1}_{10}^{\top}$ and $W=\frac{1}{20}\textbf{1}_{20}\textbf{1}_{20}^{\top}.$ The models are trained for up to 1000 epochs, which takes approximately two hours each time using a 3080 GPU. We do not incorporate dropouts in our training process.

\textbf{Dataset}.
We use two distinct datasets: CIFAR-10 and CIFAR-100. Both datasets comprise 50,000 training images and 10,000 testing images. CIFAR-10 contains images across 10 classes, while CIFAR-100 spans 100 classes. These datasets are composed of color images depicting common objects, with each image measuring 32x32 pixels with 3 color channels. Each attribute of the data
 is normalized to $[0,1]$.

\textbf{Results}.
We conducted our experiments by using the distributed mSGD with three different momentum coefficients, namely, $\alpha=0$ (corresponding to standard SGD), $\alpha=0.5$, and  $\alpha=0.9$. The experimental results, as depicted in Figures \ref{fig:test1}, illustrate some key observations:
The loss decreases to near zero across all three settings of the momentum coefficient, and the setting of $\alpha=0.9$ results in the fastest convergence of the gradient of loss to a small neighborhood around zero, outperforming the other two settings. This empirical finding is in accordance with the theoretical analysis presented in Theorems \ref{lem_mied} and \ref{dacnie}."

\section{Conclusion}
This paper presents a theoretical analysis of the last-iterate convergence properties of distributed momentum Stochastic Gradient Descent (mSGD) algorithms in non-convex settings, under the classical Robbins–Monro step-size schedule.
We establish a general analytical framework to prove both almost sure and 
$L_2$ convergence of the last iterate for a class of distributed mSGD variants, including momentum-based PSASGD, EASGD, and D-PSGD algorithms. In addition, we analyze their convergence rates.
Our results also reveal that incorporating a momentum term accelerates convergence to a neighborhood of a stationary point during the early stages of the algorithm. These findings provide important insights into the behavior and practical performance of distributed mSGD algorithms in real-world applications.
By providing the results of a classification task using the ResNet20 network, optimized with the aforementioned distributed mSGD algorithm, we observe that the experimental outcomes align well with our theoretical predictions.
In conclusion, the theoretical results presented in this work make a significant contribution to the field of distributed stochastic optimization, particularly in scenarios where computational efficiency and data privacy are critical concerns.

 \bibliographystyle{IEEEtran}
 \bibliography{ref}

\section{Biography Section}
\begin{IEEEbiography}[{\includegraphics[width=1in,height=1.25in,clip,keepaspectratio]{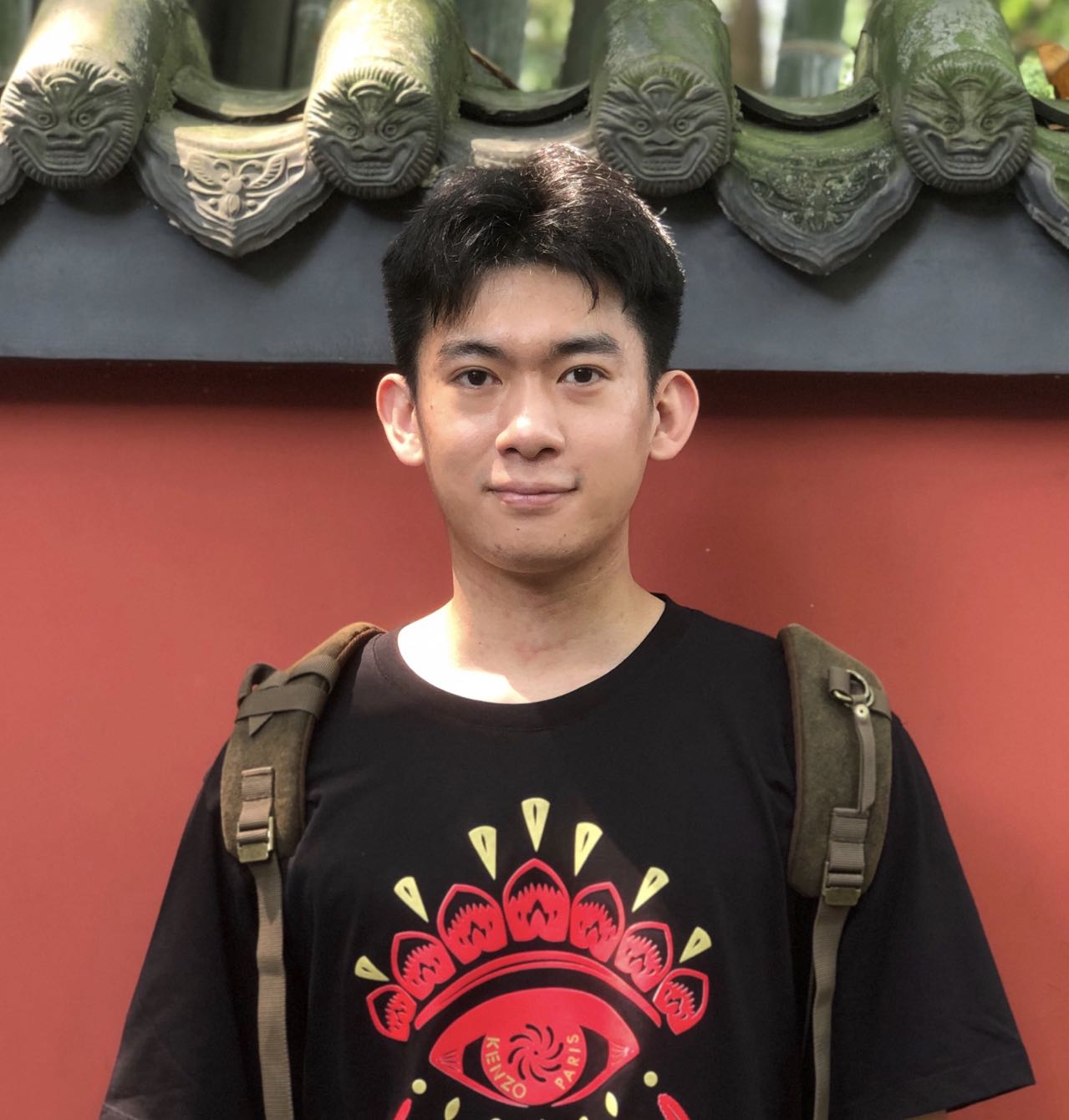}}]{Difei~Cheng} 
(M'22)received B.S. degree in mathematics and applied mathematics from Shandong University, Jinan, China, in 2017 and the Ph.D. degree in applied mathematics at the Institute of Applied Mathematics, Academy of Mathematics and Systems Science, Chinese Academy of Sciences, Beijing, China, in 2022. He is currently a Post-doctoral Research Fellow with Institute of Automation, Chinese Academy of Sciences, Beijing, China. His current research interests include clustering, stochastic optimization and deep learning.
\end{IEEEbiography} 
\begin{IEEEbiography}[{\includegraphics[width=1in,height=1.25in,clip,keepaspectratio]{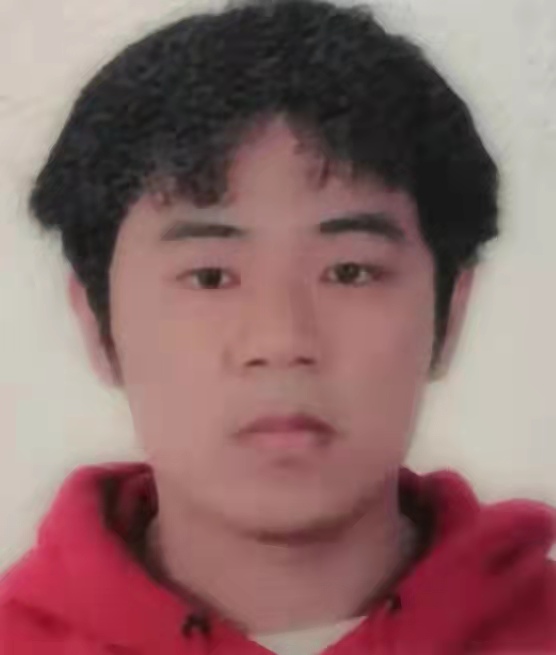}}]{Ruinan~Jin}
Ruinan Jin received the B.S.  degree in Schiffsmotor from Wuhan University of Technology, China, in 2017 and the Ph.D. degree at the Institute of Systems Science, Academy of Mathematics and Systems Science, Chinese Academy of Sciences, Beijing, China, in 2022. He is currently a Post-doctoral  Research  Fellow with School of Data Science, The Chinese University of Hong Kong, Shenzhen, China. His current research interests include  stochastic optimization, unsupervised feature learning, deep learning theory, causal discovery.

\end{IEEEbiography} 
\begin{IEEEbiography}[{\includegraphics[width=1in,height=1.25in,clip,keepaspectratio]{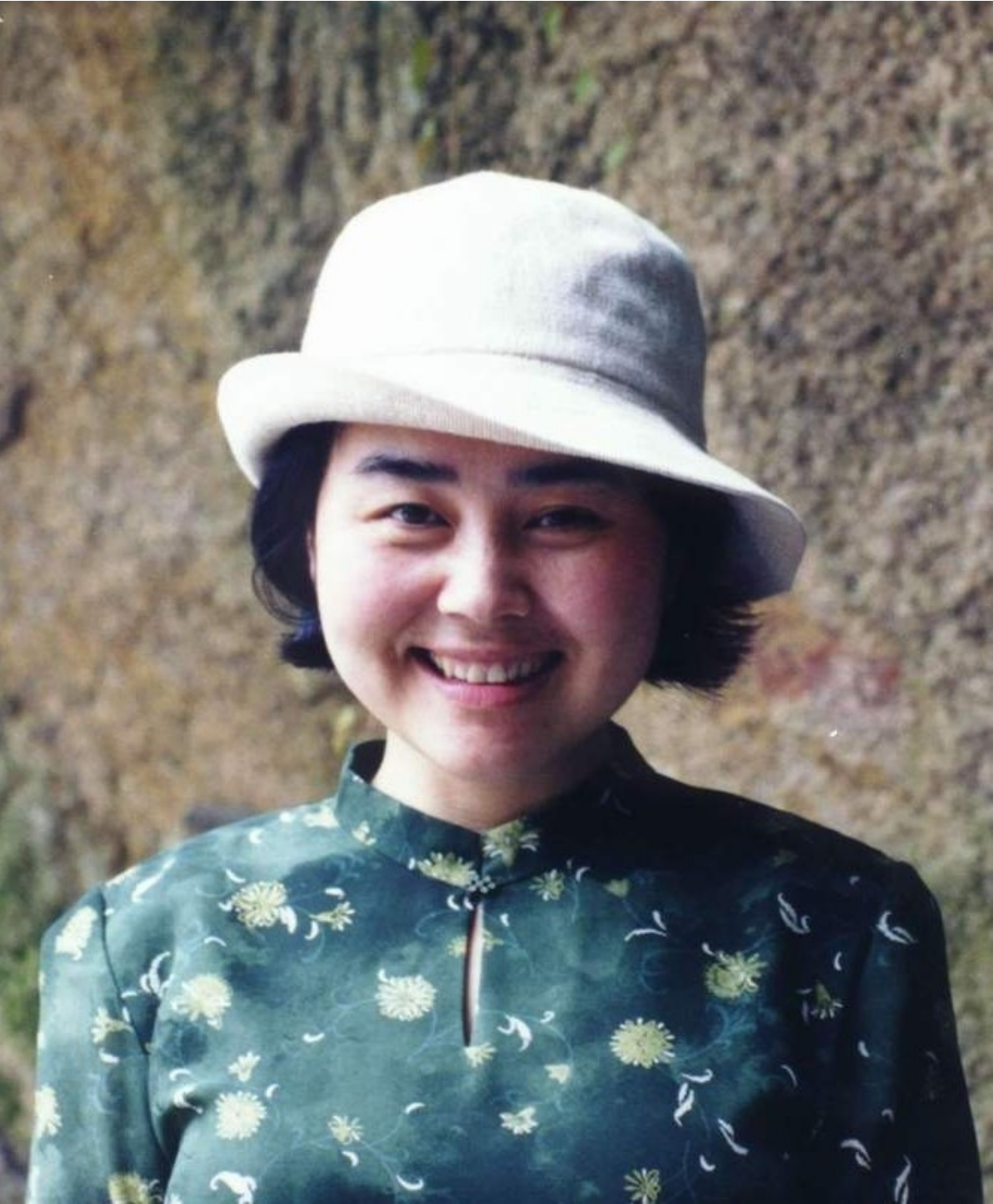}}]{Hong~Qiao}
(F'18) received the B.Eng. degree in hydraulics and control and the M.Eng. degree in robotics and automation from 
Xi'an Jiaotong University, Xi'an, China, and the Ph.D. degree in robotics from De Montfort University, Leicester, 
UK, in 1986, 1989, and 1995, respectively.
She was a University Research Fellow in De Montfort University, UK, a Research/Assistant Professor with City University 
of Hong Kong, Hong Kong, and a Lecturer with University of Manchester, Manchester, UK, from 1995 to 2004. 
She is currently a Professor with the State Key Lab of Management and Control for Complex Systems, Institute of Automation, 
Chinese Academy of Sciences, Beijing, China. 
Her current research interests include robotics, artificial intelligence, and machine learning.

 Dr. Qiao is currently a member of the Administrative Committee of
 the IEEE Robotics and Automation Society. She is the Editor-in-Chief of
 Assembly Automation, an Associate Editor of the IEEE TRANSACTIONSON
 CYBERNETICS, the IEEE TRANSACTION ON AUTOMATION SCIENCE AND ENGINEERING, the IEEE TRANSACTIONSON MECHATRONICS, the IEEE TRANSACTIONS ON COGNITIVE AND DEVELOPMENTAL SYSTEMS and the IEEE TRANSACTIONS ON NEURAL NETWORKS AND LEARNING SYSTEMS.
\end{IEEEbiography} 
\begin{IEEEbiography}[{\includegraphics[width=1in,height=1.25in,clip,keepaspectratio]{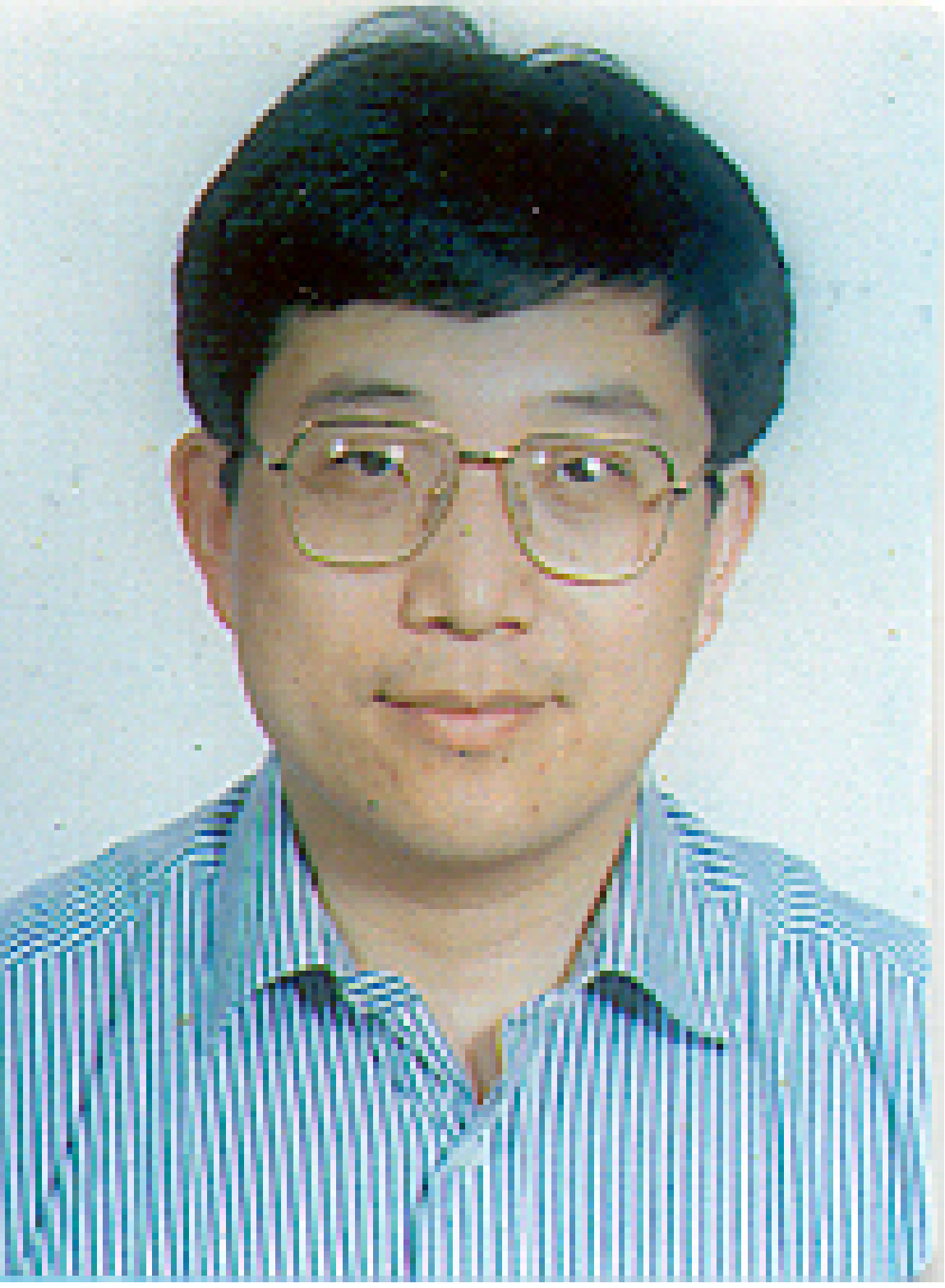}}]{Bo~Zhang}
 received B.S. degree in mathematics from Shandong University, Jinan, China, M.S. degree in mathematics from Xi'an Jiaotong University, Xi'an, China, and Ph.D. degree in applied mathematics from the University of Strathclyde, Glasgow, UK, in 1983, 1985, and 1992, respectively. After being a postdoc at Keele University, UK and a Research Fellow at Brunel University, UK, from 1992 to 1997, he joined Coventry University, Coventry, UK, in 1997, as a Senior Lecturer, where he was promoted to Reader in Applied Mathematics in 2000 and to Professor of Applied Mathematics in 2003. He is currently a Professor with Institute of Applied Mathematics, Academy of Mathematics and Systems Science, Chinese Academy of Sciences, Beijing, China. His current research interests include direct and inverse scattering problems, radar and sonar imaging, machine learning, and data mining.
\end{IEEEbiography}

\vfill

\onecolumn

\appendix
\section{Appendix}
\subsection{Useful lemmas}
\begin{lem}\label{lem_ggg}(Lemma 1.2.3,\cite{2004Introductory})
	Suppose $f(x)\in C^{1}\ \ (x\in \mathbb{R}^{N})$ with gradient satisfying the following Lipschitz condition
	\begin{equation}\nonumber\begin{aligned}
			\|\nabla  f(x)-\nabla  f(y)\|\le c\|x-y\|,
		\end{aligned}
	\end{equation}then for any $x, y\in \mathbb{R}^{N}$,  it holds that
	\begin{equation}\nonumber\begin{aligned}
			f(x)\le f(y)+\nabla  f(y)^{\top}(x-y)+\frac{c}{2}\|x-y\|^{2}.
		\end{aligned}
	\end{equation}
\end{lem}




\begin{lem}\label{lemiuy8}(Lemma 10,\cite{jin2022convergence})
	Under the same conditions as Lemma~\ref{lem_ggg}, for any  $ \ x_{0}\in \ \mathbb{R}^{N}$, it holds that
	\begin{equation}\nonumber\begin{aligned}
			\big\|\nabla  f(x_{0})\big\|^{2}\le {2c}\big(f(x_{0})-f^{*}\big), 
		\end{aligned}
	\end{equation}where $f^{*}=\inf_{x\in \ \mathbb{R}^{N}}f(x)$
\end{lem}
\begin{lem}\label{23lem_mid1} (Lemma B.6 in \cite{jin2022revisit})
	If $0<\mu<1$ and $0<\sigma<1$ ($\sigma\neq \mu$) are two constant, then exists $k_{1}>0, \ k_{2}>0$, for any positive sequence $\{\psi_{n}^{(i)}\}$, it holds that

 \begin{equation}\nonumber\begin{aligned}
			k_{1}\sum_{i=1}^{n}\kappa^{n-i}\psi_{i}\le \sum_{k=1}^{n}\mu^{n-k}\sum_{i=1}^{k}\sigma^{k-i}\psi_{i}\le k_{2}\sum_{i=1}^{n}\kappa^{n-i}\psi_{i},
	\end{aligned}\end{equation}
     
  where $\kappa=\max\{\mu,\sigma\}$ and $\omega_{0}=\log_{\kappa}\min\{\mu,\sigma\}$.
\end{lem}	
\begin{lem}\label{lem_small} 
If there exists a sequence of positive numbers $\{x_n\}_{n=1}^{\infty}$ such that $\sum_{n=1}^{\infty} x_n < \infty$, then for any $n > 0$, there exists a constant $k_n > 0$, uniform in $n$, such that for any $s$, it holds that $\sum_{k=s}^{n} x_k < k_n x_s$.
\end{lem}

\begin{lem} \label{lem_summation_MDS}\cite{wang2019almost}
	Suppose that $\{X_{n}\}\in \mathbb{R}^{N}$ is a $\mathcal{L}_2$ martingale difference sequence, and $(X_{n},\mathcal{ F}_{n})$ is an adaptive process. Then it holds that $\sum_{k=0}^{\infty}X_k<+\infty \ a.s.,$ if  $\sum_{n=1}^{\infty}\Expect(\|X_{n}\|^{2})<+\infty$ or $\sum_{n=1}^{\infty}\Expect\big(\|X_{n}\|^{2}\big|\mathscr{F}_{n-1}\big)<+\infty.$
\end{lem}		

\begin{lem}\label{lem_summation}(Lemma 6,\cite{jin2022convergence})
	Suppose that $\{X_n\}\in \mathbb{R}^{N}$ is a non-negative sequence of random variables, then it holds that $\sum_{n=0}^{\infty}X_n<+\infty \ a.s.,$ if   
	$\sum_{n=0}^{\infty}\Expect\big(X_n\big)<+\infty.$   
\end{lem}
\subsection{Proof of Theorem \ref{lem_mied}}

\begin{proof}
	First, due to $0<\alpha<1,$ we can always find a positive constant $\alpha_{0},$ making $\alpha_{1}:=\alpha^{2}+\alpha_{0}<1.$ Then based on the first equation of E.q. \eqref{123combinemomentum}, we have
	\begin{equation}\label{-0}\begin{aligned}
			\|v_{n}\|^{2}\le \alpha_{1} \|v_{n-1}\|^{2}+\epsilon^{2}_{n}\Big(1+\frac{1}{{\alpha_{0}}}\Big)\cdot\|G_{n,\xi_{n}}\|^{2}.
		\end{aligned} 
	\end{equation} Then we take the mathematical expectation on the both side of E.q. \eqref{-0}, acquiring

	\begin{equation}\nonumber\begin{aligned}
			\Expect\|v_{n}\|^{2}&\le \alpha_{1}\cdot \Expect\|v_{n-1}\|^{2}+\epsilon^{2}_{n}\Big(1+\frac{1}{\alpha_{0}}\Big)\cdot\Expect\|G_{n,\xi_{n}}-G_{n}\|^{2}+\epsilon^{2}_{n}\Big(1+\frac{1}{\alpha_{0}}\Big)\cdot\Expect\|G_{n}\|^{2}\\&\le  \alpha_{1}\cdot \Expect\|v_{n-1}\|^{2}+\epsilon^{2}_{n}\sigma_{0}^{2}\Big(1+\frac{1}{\alpha_{0}}\Big)+\epsilon^{2}_{n}\Big(1+\frac{1}{\alpha_{0}}\Big)\cdot\Expect\|G_{n}\|^{2}.
		\end{aligned} 
	\end{equation}
	Iterating above inequity, we acquire
	
	\begin{equation}\label{msgd_03}
		\begin{aligned}
			\Expect\|v_{n}\|^{2}\le \alpha_{1}^{n}\cdot\|v_{0}\|^{2}+\epsilon_{n}^{2}\Big(1+\frac{1}{\alpha_{0}}\Big)\cdot\sum_{s=1}^{n}(\sigma_{0}^{2}+\Expect\|G_{s}\|^{2})\cdot\alpha_{1}^{n-s}
		\end{aligned}
	\end{equation}
		Next, we iterate the second equation of E.q. \eqref{123combinemomentum} to attain 

	\begin{equation}\nonumber
		\begin{aligned}
			X_{n+1}&=W_{n}X_{n}-W_{n}v_{n}\\
			&=W_{n}\big(W_{n-1}X_{n-1}-W_{n-1}v_{n-1}\big)-W_{n}v_{n}\\&=W_{n}W_{n-1}(W_{n-2}X_{n-2}-W_{n-2}v_{n-2})\\
			&-W_{n}W_{n-1}v_{n-1}-W_{n}v_{n}\\&=...=\bigg(\prod_{s=1}^{n}W_{s}\bigg)\cdot X_{1}-\sum_{t=1}^{n}\Bigg(\bigg(\prod_{s=t}^{n}W_{s}\bigg)\cdot v_{t}\Bigg).
		\end{aligned}
	\end{equation}
 We left multiply both sides of the above equation by the vector $e_i^{\top}:=(-1/m,-1/m,...,1-1/m,...,-1/m,...,-1/m)$ ( the $i$-th entry is  $1-1/m$, and others are $1/m$) to obtain

	\begin{equation}\label{msgd_001}
		\begin{aligned}
			e_{i}^{\top}X_{n+1}=e_{i}^{\top}\bigg(\prod_{s=1}^{n}W_{s}\bigg)\cdot X_{1}-\sum_{t=1}^{n}\Bigg(e_{i}^{\top}\bigg(\prod_{s=t}^{n}W_{s}\bigg)\cdot v_{t}\Bigg).
		\end{aligned}
	\end{equation}
 When $s \mod k = 0$, according to assumption\ref{ass_g1} 2), we can find an orthogonal matrix $Q$ such that \(Q^{\top} W_{s}Q = \text{diag}\{1, \lambda_{2}, \lambda_{3}, \ldots, \lambda_{m}\},\) where $\lambda_{0}:=\max_{2\le j\le m}\{|\lambda_{j}|\}<1.$ When $n \mod k \neq 0,$ we always have $W_s = \textbf{I}_{m} = QQ^{\top}.$ We assign $W_{t,n}:=\prod_{s=t}^{n}W_{s}.$ Then we can get that
 
 		\begin{equation}\label{msgd_01}\begin{aligned}
 			W_{t,n}&=\prod_{s\in [t,n],\ s\mod k=0}W_{s}\\
 			&=\prod_{s\in [t,n],\ s\mod k=0}(Q\cdot\text{diag}\{1,\lambda_{2}, \lambda_{3}, \ldots, \lambda_{m}\}\cdot Q^{\top})\\&=Q\cdot\text{diag}\bigg\{1,\prod_{s\in [t,n],\ s\mod k=0}\lambda_{2}, \\
 			&\prod_{s\in [t,n],\ s\mod k=0}\lambda_{3}, \ldots, \prod_{s\in [t,n],\ s\mod k=0}\lambda_{m}\bigg\}\cdot Q^{\top}.
 	\end{aligned}\end{equation}
  We can conclude that $\forall\ j\in[2,m],$ there is
	\begin{equation}\label{msgd_12}\begin{aligned}
			\Bigg|\prod_{s\in [t,n],\ s\mod k=0}\lambda_{j}\Bigg|\le \lambda_{0}^{c(t,n)},\end{aligned}\end{equation} where $c(t, n)$ represents the total number of integers divisible by $k$ between $t$ and $n$. It is easy to prove that $$\bigg\lfloor\frac{n-t}{k}\bigg\rfloor\le c(t,n)\le \bigg\lfloor\frac{n-t}{k}\bigg\rfloor+1.$$
	Then based on E.q. \eqref{msgd_01}, we can derive the following expression:
	\begin{equation}\label{msgd_02}\begin{aligned}
			&\|e_{i}^{\top}X_{n+1}\|^{2}\le 2\|e_{i}^{\top}W_{1,n}\|^{2}\cdot \|X_{1}\|^{2}+2\bigg(\sum_{t=1}^{n}\|e_{i}^{\top}W_{t,n}\|\cdot \|v_{t}\|\bigg)^{2}.
	\end{aligned}\end{equation}
	For any \(t\) and \(n\), since the matrix \(W_{t,n}\) is a real symmetric matrix, its eigenspaces are orthogonal to each other. We know that \((1,1,...,1)^\top\) is obviously an eigenvector corresponding to the eigenvalue 1, and according to Assumption \ref{ass_g1} 2), we know that the dimension of the eigenspace corresponding to the eigenvalue 1 can only be 1. Therefore, the eigenspace corresponding to the eigenvalue 1 is completely spanned by the vector \((1,1,...,1)^\top\). On the other hand, since \(e_{i}^\top (1,1,...,1)^\top = 0\), we know that \(e_{i}\) must belong to the direct sum of the eigenspaces of \(W_{t,n}\) other than the one corresponding to eigenvalue 1. Hence, there exists an orthogonal decomposition $$e_{i} = r_2 e_{2,i} + r_3 e_{3,i} + ... + r_m e_{m,i},$$ where each \(e_{s,i}\) (\(2 \le s \le m\)) is a unit vector and at the same time an eigenvector of the matrix \(W_{t,n}\) corresponding to an eigenspace not associated with the eigenvalue 1. Therefore, we can obtain
	\begin{equation}
		\begin{aligned}
			\|e_{i}^\top W_{t,n}\| &= \bigg\|\sum_{s=2}^{m} r_s e_{s,i}^\top W_{t,n}\bigg\| = \bigg(\sum_{s=2}^{m} |r_s^2|\bigg)^{\frac{1}{2}} \cdot \lambda_0^{c(t,n)} \\
			&= \|e_{i}\|\lambda_0^{c(t,n)}.
		\end{aligned}		
	\end{equation}
Substitute above inequity into E.q. \eqref{msgd_02}, getting
	
	    \begin{equation}\nonumber
			\begin{aligned}
				\|e_{i}^{\top}X_{n+1}\|^{2}&\le 2\|e_{i}\|^{2}\cdot\|X_{1}\|^{2}\cdot\lambda_{0}^{2c(1,n)}+2\|e_{i}\|^{2}\bigg(\sum_{t=1}^{n} \lambda_0^{c(t,n)}\cdot \|v_{t}\|\bigg)^{2}\\&\le 2\|e_{i}\|^{2}\cdot\|X_{1}\|^{2}\cdot\lambda_{0}^{c(1,n)}+2\|e_{i}\|^{2}\lambda(n,k)\sum_{t=1}^{n} \lambda_{0}^{c(t,n)}\cdot \|v_{t}\|^{2},
				\end{aligned}
			\end{equation}	
	where $\lambda(n,k)=\sum_{t=1}^{n}{\lambda_{0}^{c(t,n)}}$. We take the mathematical expectation, resulting
	\begin{equation}\label{msgd_04}\begin{aligned}
			&\Expect\|e_{i}^{\top}X_{n+1}\|^{2}\le 2\|e_{i}\|^{2}\cdot\|X_{1}\|^{2}\cdot\lambda_{0}^{c(1,n)}+2\|e_{i}\|^{2}\cdot\lambda(n,k)\sum_{t=1}^{n} \lambda_{0}^{c(t,n)}\cdot \Expect\|v_{t}\|^{2}.
	\end{aligned}\end{equation}
	We substitute \eqref{msgd_03} into \eqref{msgd_04}, getting
	
			\begin{equation}\nonumber\begin{aligned}
				&\Expect\|e_{i}^{\top}X_{n+1}\|^{2}\le 2\|e_{i}\|^{2}\cdot\|X_{1}\|^{2}\cdot\lambda_{0}^{c(1,n)}+2\|e_{i}\|^{2}\cdot\lambda(n,k)\\
				&\sum_{t=1}^{n} \lambda_{0}^{c(t,n)}\cdot \bigg(\alpha_{1}^{t}\|v_{0}\|^{2}+\epsilon_{t}^{2}\Big(1+\frac{1}{\alpha_{0}}\Big)\sum_{s=1}^{t}(\sigma_{0}^{2}+\Expect\|G_{s}\|^{2})\cdot\alpha_{1}^{t-s}\bigg)\\&=2\|e_{i}\|^{2}\cdot\|X_{1}\|^{2}\cdot\lambda_{0}^{c(1,n)}+2\|e_{i}\|^{2}\lambda(n,k)\sum_{t=1}^{n} \lambda_{0}^{c(t,n)}\cdot\alpha_{1}^{t}\cdot\|v_{0}\|^{2}\\&+2\|e_{i}\|^{2}\lambda(n,k)\sum_{t=1}^{n}\lambda_{0}^{c(t,n)}\cdot\epsilon_{t}^{2}\Big(1+\frac{1}{\alpha_{0}}\Big)\sum_{s=1}^{t}(\sigma_{0}^{2}+M^{2})\alpha_{1}^{t-s}.
		\end{aligned}\end{equation}
	In the above inequality, by substituting the estimate for $c(t,n)$ from \eqref{msgd_12} and simplifying, we can obtain
	\begin{equation}\label{mSGD_13}\begin{aligned}
			&\Expect\|e_{i}^{\top}X_{n+1}\|^{2}=O\bigg(\sum_{t=1}^{n}\max\{\lambda_{0}^{\frac{1}{k}},\alpha_{1}\}^{n-t}\cdot\epsilon_{t}^{2}\bigg)\rightarrow 0.
	\end{aligned}\end{equation}
	We recall E.q. \eqref{123combinemomentum} as follows
	\begin{equation}\nonumber\begin{aligned}
			&v_{n}=\alpha v_{n-1}+\epsilon_{n}G(X_{n},\xi_{n}), \\&X_{n+1}= W_{n}\big(X_{n}-v_{n}\big).
	\end{aligned} \end{equation}Then we multiply $u^{\top}=(1/m,1/m,...,1/m)$ on the both sides of the above equalities to obtain 
	\begin{equation}\nonumber\begin{aligned}
			&u^{\top}v_{n}=\alpha u^{\top}v_{n-1}+\epsilon_{n}u^{\top}G(X_{n},\xi_{n}), \\& u^{\top}X_{n+1}=u^{\top} W_{n}\big(X_{n}-v_{n}\big).
	\end{aligned}\end{equation} Since $ W_{n}$ is a doubly stochastic matrix,    $u^{\top} W_{n}=u^{\top}$, Furthermore, it holds that 
	\begin{equation}\nonumber\begin{aligned}
			&u^{\top}v_{n}=\alpha u^{\top}v_{n-1}+\epsilon_{n}u^{\top}G(X_{n},\xi_{n}), \\& u^{\top}X_{n+1}=u^{\top}X_{n}-u^{\top}v_{n}.
	\end{aligned}\end{equation}
 Denote $\textbf{I}_{m}=[1,1,\cdots,1]_{1\times m}$ and $\otimes$ is Kronecker Product. We derive $g(u^{\top}X_{n+1})-g( u^{\top}X_{n})$ to obtain
	\begin{equation}\label{podjc}\begin{aligned}
			g(u^{\top}X_{n+1})-g( u^{\top}X_{n})&=-(u^{\top}G(\textbf{I}_{m}\otimes u^{\top}X_{n}))^{\top}(u^{\top}v_{n})+(u^{\top}G(\textbf{I}_{m}\otimes u^{\top}X_{n})-u^{\top}G(\textbf{I}_{m}\otimes u^{\top}X_{\zeta_{n}}))^{\top}(u^{\top}v_{n})\\&\le-(u^{\top}G(\textbf{I}_{m}\otimes u^{\top}X_{n}))^{\top}(u^{\top}v_{n})+L\|u^{\top}v_{n}\|^{2},
	\end{aligned}\end{equation} 
	where $u^{\top}X_{\zeta_{n}}$ is a value between $u^{\top}X_{n}$ and $u^{\top}X_{n+1}$. Next we focus on the term $(u^{\top}G(\textbf{I}_{m}\otimes u^{\top}X_{n}))^{\top}(u^{\top}v_{n})$. We derive that
 
 \begin{equation}\label{ionjuhf}\begin{aligned}
			&(u^{\top}G(\textbf{I}_{m}\otimes u^{\top}X_{n}))^{\top}(u^{\top}v_{n})=(u^{\top}G(\textbf{I}_{m}\otimes u^{\top}X_{n}))^{\top}(\alpha u^{\top}v_{n-1}+\epsilon_{n}u^{\top}G(X_{n},\xi_{n}))\\&=\alpha(u^{\top}G(\textbf{I}_{m}\otimes u^{\top}X_{n}))^{\top}u^{\top}v_{n-1}+\epsilon_{n}(u^{\top}G(\textbf{I}_{m}\otimes u^{\top}X_{n}))^{\top}u^{\top}G(X_{n},\xi_{n}))\\&\ge \alpha u^{\top}G(\textbf{I}_{m}\otimes u^{\top}X_{n-1})^{\top}(u^{\top}v_{n-1})-L\|u^{\top}v_{n-1}\|^{2}+\epsilon_{n}(\textbf{I}_{m}\otimes u^{\top}G(X_{n}))^{\top}u^{\top}G(X_{n},\xi_{n})).
	\end{aligned}\end{equation}   
 It follows from E.q. \eqref{ionjuhf} that
	\begin{equation}\label{imksodjw}\begin{aligned}
			(u^{\top}G(\textbf{I}_{m}\otimes u^{\top}X_{n}))^{\top}(u^{\top}v_{n})\ge-L\sum_{s=0}^{n-1}\alpha^{n-s-1}\|u^{\top}v_{s}\|^{2}+\sum_{s=1}^{n}\alpha^{n-s}\epsilon_{s}(\textbf{I}_{m}\otimes u^{\top}G(X_{s}))^{\top}u^{\top}G(X_{s},\xi_{s})).
	\end{aligned}\end{equation} Substituting E.q. \eqref{imksodjw} into E.q. \eqref{podjc} leads to

	\begin{equation}\label{cxerqwe}\begin{aligned}
			g(u^{\top}X_{n+1})-g( u^{\top}X_{n})\le {L}\sum_{s=1}^{n}\alpha^{n-s}\|u^{\top}v_{s}\|^{2}-\sum_{s=1}^{n}\alpha^{n-s}\epsilon_{s}(\textbf{I}_{m}\otimes u^{\top}G(X_{s}))^{\top}u^{\top}G(X_{s},\xi_{s}))+L\|u^{\top}v_{0}\|^{2}\cdot\alpha^{n}.
	\end{aligned}\end{equation}
	 Then we consider the term $(u^{\top}G(\textbf{I}_{m}\otimes u^{\top}X_{s}))^{\top}u^{\top}G(X_{s},\xi_{s}))$ to have 
	 
   \begin{equation}\label{fderw}\begin{aligned}
	 			&-(u^{\top}G(\textbf{I}_{m}\otimes u^{\top}X_{s}))^{\top}u^{\top}G(X_{s},\xi_{s}))\\&=-u^{\top}G(\textbf{I}_{m}\otimes u^{\top}X_{s})^{\top}\big(u^{\top}G(X_{s},\xi_{s})-u^{\top}G(X_{s})\big)-\|(u^{\top}G(\textbf{I}_{m}\otimes u^{\top}X_{s}))^{\top}\|^{2}\\&+u^{\top}G(\textbf{I}_{m}\otimes u^{\top}X_{s})^{\top}(u^{\top}G(\textbf{I}_{m}\otimes u^{\top}X_{s})-u^{\top}G(X_{s}))\\&\le-\frac{1}{2}\big\|u^{\top}G(\textbf{I}_{m}\otimes u^{\top}X_{s})\big\|^{2}+2L\sum_{i=1}^{m}\|x_{s}^{(i)}-u^{\top}X_{s}\|^{2}-u^{\top}G(\textbf{I}_{m}\otimes u^{\top}X_{s})^{\top}\big(u^{\top}G(X_{s},\xi_{s})-u^{\top}G(X_{s})\big).
	 	\end{aligned}\end{equation} 
	 Denote $\beta_{s}:=2L\sum_{i=1}^{m}\|x_{s}^{(i)}-u^{\top}X_{s}\|^{2}$, then substituting E.q. \eqref{fderw} into E.q. \eqref{cxerqwe} yields
	
		\begin{equation}\label{vcdf}\begin{aligned}
				&g(u^{\top}X_{n+1})-g( u^{\top}X_{n})\le {L}\sum_{s=1}^{n}\alpha^{n-s}\|u^{\top}v_{s}\|^{2}-\frac{1}{2}\sum_{s=1}^{n}\alpha^{n-s}\epsilon_{s}\big\|u^{\top}G(\textbf{I}_{m}\otimes u^{\top}X_{s})\big\|^{2}+\sum_{s=1}^{n}\alpha^{n-s}\epsilon_{s}\beta_{s}\\
    &-\sum_{s=1}^{n}\alpha^{n-s}\epsilon_{s}u^{\top}G(\textbf{I}_{m}\otimes u^{\top}X_{s})^{\top}\cdot\big(u^{\top}G(X_{s},\xi_{s})-u^{\top}G(X_{s})\big)+L\|u^{\top}v_{0}\|^{2}\cdot\alpha^{n}.
		\end{aligned}\end{equation}
		On the other hand, we have
	
		\begin{equation}\label{dqsc}\begin{aligned}
				\|u^{\top}v_{n}\|^{2}&=\|\alpha u^{\top}v_{n-1}+\epsilon_{n}u^{\top}G(X_{n},\xi_{n})\|^{2}\\
				&=\alpha^{2}\|u^{\top}v_{n-1}\|^{2}+2\alpha\epsilon_{n}(u^{\top}v_{n-1})^{\top}u^{\top}G(X_{n},\xi_{n})+\epsilon_{n}^{2}\|u^{\top}G(X_{n},\xi_{n})\|^{2}\\
				&=\alpha^{2}\|u^{\top}v_{n-1}\|^{2}+2\alpha\epsilon_{n}(u^{\top}v_{n-1})^{\top}u^{\top}G(X_{n})+\epsilon_{n}^{2}\|u^{\top}G(X_{n},\xi_{n})\|^{2}+\gamma_{n},
		\end{aligned}\end{equation}
	 where $\gamma_{n}=2\alpha\epsilon_{n}v_{n-1}^{\top}u^{\top}(G(X_{n},\xi_{n})-G(X_{n}))$. Then we calculate $2\epsilon_{n}(E.q. \eqref{podjc}-E.q. \eqref{ionjuhf})+E.q. \eqref{dqsc}$ to obtain

		\begin{equation}\label{daceio}\begin{aligned}
			&2\epsilon_{n+1}g(u^{\top}X_{n+1})-2\epsilon_{n}g(u^{\top}X_{n})\\
   &\le \alpha^{2}\|u^{\top}v_{n-1}\|^{2}-\|u^{\top}v_{n}\|^{2}+2\epsilon_{n}L\|u^{\top}v_{n}\|^{2}+\hat{\gamma}_{n}+\epsilon_{n}^{2}\|u^{\top}G(X_{n},\xi_{n})\|^{2}-2\epsilon_{n}^{2}(u^{\top}G(\textbf{I}_{m}\otimes u^{\top}X_{n}))^{\top}u^{\top}G(X_{n},\xi_{n})\\&\le\alpha^{2}\|u^{\top}v_{n-1}\|^{2}-\|u^{\top}v_{n}\|^{2}+2\epsilon_{n}L(\|u^{\top}v_{n}\|^{2}+\|u^{\top}v_{n-1}\|^{2})+\epsilon_{n}^{2}\|u^{\top}G(X_{n},\xi_{n})\|^{2}\\&+\hat{\gamma}_{n}-\epsilon_{n}^{2}\|u^{\top}G(\textbf{I}_{m}\otimes u^{\top}X_{n})\|^{2}+2\epsilon_{n}^{2}\beta_{n},
	\end{aligned}\end{equation}
where 

		\begin{equation}\nonumber
		\hat{\gamma}_{n}:=\gamma_{n}+2\epsilon_{n}^{2}(u^{\top}G(\textbf{I}_{m}\otimes u^{\top}X_{n}))^{\top}\big(u^{\top}G(X_{n},\xi_{n})-u^{\top}G(X_{n})\big).
	\end{equation}
We make the mathematical expectation of E.q. \eqref{vcdf} to obtain
 
 		\begin{equation}\nonumber\begin{aligned}
 			&\Expect\big(g(u^{\top}X_{n+1})\big)-\Expect\big(g(u^{\top}X_{n})\big)\le \hat{L}\sum_{s=1}^{n}\alpha^{n-s}\Expect\|u^{\top}v_{s}\|^{2}-\frac{1}{2}\sum_{s=1}^{n}\alpha^{n-s}\epsilon_{s}\Expect\big\|u^{\top}G(\textbf{I}_{m}\otimes u^{\top}X_{s})\big\|^{2}\\			&+L\|u^{\top}v_{0}\|\cdot\alpha^{n}+\sum_{s=1}^{n}\alpha^{n-s}\epsilon_{s}\beta_{s}.
 	\end{aligned}\end{equation}
  Making a summation of the above inequality leads to 
 
 	\begin{equation}\label{dsaucheiuas}\begin{aligned}
 			&\Expect\big(g(u^{\top}X_{n+1})\big)-\Expect\big(g(u^{\top}X_{1})\big)\le \frac{L}{1-\alpha}\sum_{s=1}^{n}\Expect\|u^{\top}v_{s}\|^{2}-\frac{1}{2}\sum_{s=1}^{n}\epsilon_{s}\Expect\big\|u^{\top}G(\textbf{I}_{m}\otimes u^{\top}X_{s})\big\|^{2}+\frac{L\|u^{\top}v_{0}\|^{2}}{1-\alpha}+\hat{\beta}_{n},
 	\end{aligned}\end{equation} where $\hat{\beta}_{n}=\sum_{t=1}^{n}\sum_{s=1}^{\top}\alpha^{t-s}\epsilon_{s}\beta_{s}$. We perform the same operations on E.q. \eqref{daceio} to obtain

  \begin{equation}\label{ocneownvjfis'}\begin{aligned}
	 			2\epsilon_{n+1}\Expect\big(g(u^{\top}X_{n+1})\big)-2\epsilon_{1}\Expect\big(g(u^{\top}X_{1})\big)\le -\sum_{s=1}^{n}(1-\alpha^{2})\Expect\|u^{\top}v_{s}\|^{2}+\sum_{s=1}^{n}\epsilon_{s}^{2}\Expect\|u^{\top}G(X_{s},\xi_{s})\|^{2}+2\sum_{s=1}^{n}\epsilon_{s}^{2}\beta_{s}.
	 	\end{aligned}\end{equation} For the term $\sum_{s=1}^{n}\epsilon_{s}^{2}\Expect\|u^{\top}G(X_{s},\xi_{s})\|^{2},$ we have
	
			\begin{equation}
			\nonumber\begin{aligned}
				&\sum_{s=1}^{n}\epsilon_{s}^{2}\Expect\|u^{\top}G(X_{s},\xi_{s})\|^{2}\le 2\sum_{s=1}^{n}\epsilon_{s}^{2}\|u^{\top}G(X_{s},\xi_{s})-u^{\top}G(X_{s})\|^{2}+2\sum_{s=1}^{n}\epsilon_{s}^{2}\Expect\|u^{\top}G(X_{s})\|^{2}\\&\le 2\sum_{s=1}^{n}\epsilon_{s}^{2}\|u^{\top}G(X_{s},\xi_{s})-u^{\top}G(X_{s})\|^{2}+4\sum_{s=1}^{n}\epsilon_{s}^{2}\Expect\|u^{\top}G(X_{s})-u^{\top}G(\textbf{I}_{m}\otimes u^{\top}X_{s})\|^{2}\\&+4\sum_{s=1}^{n}\epsilon_{s}^{2}\Expect\|u^{\top}G(\textbf{I}_{m}\otimes u^{\top}X_{s})\|^{2}.
		\end{aligned}\end{equation}
	 From Assumption \ref{ass_g1} Item (4), we know that
	\begin{equation}\nonumber
		\begin{aligned}
		2\sum_{s=1}^{n}\epsilon_{s}^{2}\|u^{\top}G(X_{s},\xi_{s})-u^{\top}G(X_{s})\|^{2}+4\sum_{s=1}^{n}\epsilon_{s}^{2}\Expect\|u^{\top}G(X_{s})-u^{\top}G(\textbf{I}_{m}\otimes u^{\top}X_{s})\|^{2}\le 2(\sigma_{0}^{2}+2L\sigma_{1})\sum_{s=1}^{n}\epsilon_{s}^{2},
	\end{aligned}\end{equation} which means

		\begin{equation}\nonumber\begin{aligned}
			\sum_{s=1}^{n}\epsilon_{s}^{2}\Expect\|u^{\top}G(X_{s},\xi_{s})\|^{2}\le 2(\sigma_{0}^{2}+2L\sigma_{1})\sum_{s=1}^{n}\epsilon_{s}^{2}+4\sum_{s=1}^{n}\epsilon_{s}^{2}\Expect\|u^{\top}G(\textbf{I}_{m}\otimes u^{\top}X_{s})\|^{2}.
	\end{aligned}\end{equation}
Substitute above inequity into E.q. \eqref{ocneownvjfis'}, getting
	\begin{equation}\label{ocneownvjfis}\begin{aligned}
			&2\epsilon_{n+1}\Expect\big(g(u^{\top}X_{n+1})\big)-2\epsilon_{1}\Expect\big(g(u^{\top}X_{1})\big)\le -\sum_{s=1}^{n}(1-\alpha^{2})\Expect\|u^{\top}v_{s}\|^{2}+2(\sigma_{0}^{2}+2\sigma_{1})\sum_{s=1}^{n}\epsilon_{s}^{2}\\&+4\sum_{s=1}^{n}\epsilon_{s}^{2}\Expect\|u^{\top}G(\textbf{I}_{m}\otimes u^{\top}X_{s})\|^{2}+2\sum_{s=1}^{n}\epsilon_{s}^{2}\beta_{s}.
	\end{aligned}\end{equation}
	We calculate $\frac{1-\alpha}{L}E.q. \eqref{dsaucheiuas}+\frac{1}{1-\alpha^{2}}E.q. \eqref{ocneownvjfis}$, from Assumption \ref{ass_g1} 4) and E.q. \eqref{mSGD_13} ($\sum_{s=1}^{n}\epsilon_{s}^{2}\beta_{s}\rightarrow 0,$ $\hat{\beta}_{n}\rightarrow 0$), we can get
	
	    \begin{equation}\nonumber\begin{aligned}
				\sum_{s=1}^{+\infty}\epsilon_{s}\Expect\|\nabla g(\overline{x}_{n})\|^{2}<+\infty,\ \  \sum_{s=1}^{+\infty}\epsilon_{s}\|\nabla g(\overline{x}_{n})\|^{2}<+\infty\ \ a.s.,
		\end{aligned}\end{equation} where the second inequity is because Lemma \ref{lem_summation}. Then by using the condition $\sum_{n=1}^{+\infty}\epsilon_{n}=+\infty,$ we can immediately acquire
	 	
	 	\begin{equation}\nonumber
	 		\liminf_{n\rightarrow+\infty}\Expect\|\nabla g(\overline{x}_{n})\|^{2}=0, \ \ \liminf_{n\rightarrow+\infty}\|\nabla g(\overline{x}_{n})\|^{2}=0\ a.s.\ .
	 	\end{equation}
	 Our goal below is to prove
	
		\begin{equation}\nonumber
			\limsup_{n\rightarrow+\infty}\Expect\|\nabla g(\overline{x}_{n})\|^{2}=0, \ \ \limsup_{n\rightarrow+\infty}\|\nabla g(\overline{x}_{n})\|^{2}=0\ a.s.\ .
		\end{equation}
	We first prove $\limsup_{n\rightarrow+\infty}\|\nabla g(\overline{x}_{n})\|^{2}=0\ a.s.\ .$ We use proof by contradiction. We assume that for a certain trajectory $\{\|\nabla g(\overline{x}_{n})\|^{2}\}_{n=1}^{+\infty},$ apart from $0$, there exists another accumulation point $\hat{u}>0\ .$ Then, for a certain open interval $(o,e) \subset (0,\hat{u}),$ the sequence $\{\| \nabla g(\overline{x}_{n})\|^{2}\}_{n=1}^{+\infty}$ must cross this interval infinitely many times. We denote all the intervals that go upwards as $\{(\|\nabla g(\overline{x}_{l_n})\|^{2},\|\nabla g(\overline{x}_{r_n})\|^{2})\}_{n=1}^{+\infty}\ .$ We have
	\begin{equation}\label{pm,41}\begin{aligned}
			\sum_{n=1}^{+\infty}\sum_{i=l_{n}}^{r_{n}}\epsilon_{i}< \frac{1}{o}\sum_{n=1}^{+\infty}\sum_{i=l_{n}}^{r_{n}}\epsilon_{i}\|\nabla g(\overline{x}_{i})\|^{2}<+\infty\ .
	\end{aligned}\end{equation} On the other hand, due to $\|\nabla g(\overline{x}_{r_{n}})\|^{2}>e$ and $\|\nabla g(\overline{x}_{l_{n}})\|^{2}<e,$ we know there is a $\tilde{p}_{0}>0,$ such that $\|\theta_{r_{n}}-\theta_{l_{n}}\|>\tilde{p}_{0}\ .$ Then we get
	\begin{equation}\nonumber\begin{aligned}
			\tilde{p}_{0}<\|\theta_{r_{n}}-\theta_{l_{n}}\|=\zeta_{n}+k_{0}\sum_{i=l_{n}}^{r_{n}}\epsilon_{i},
	\end{aligned}\end{equation} where $\zeta_{n}\rightarrow \ 0\ .$ We get
	$$\liminf_{n\rightarrow+\infty}\sum_{i=l_{n}}^{r_{n}}\epsilon_{i}>\frac{\tilde{p}_{0}}{2k_{0}}>0,$$ which conclude
	\begin{equation}\label{dfsafsad}\begin{aligned}
			\sum_{n=1}^{+\infty}\sum_{i=l_{n}}^{r_{n}}\epsilon_{i}=+\infty\ .
	\end{aligned}\end{equation} Now we have a contradiction between E.q. \eqref{dfsafsad} and E.q. \eqref{pm,41}, which implies that our assumption is false. Therefore, we obtain $\limsup_{n\rightarrow+\infty}\|\nabla g(\overline{x}_{n})\|^{2}=0\ a.s.,$ that is $\lim_{n\rightarrow+\infty}\|\nabla g(\overline{x}_{n})\|^{2}=0\ a.s.\ .$ Using the same technique, we can obtain convergence in the mean square sense, i.e., $\lim_{n\rightarrow+\infty}\Expect\|\nabla g(\theta_{n})\|^{2}=0$ from the inequity $\sum_{s=1}^{+\infty}\epsilon_{s}\Expect\|\nabla g(\overline{x}_{n})\|^{2}<+\infty\ .$
\end{proof}
\subsection{Proof of Theorem \ref{coro_41}}
\begin{proof}
We define $z_{n}=\frac{u^{\top}X_{n}-\alpha u^{\top}X_{n-1}}{1-\alpha}.$ We can obtain $\forall\ \theta_{0}\in\mathbb{R}^{d}$ which satisfies $\|z_{n}-\theta_{0}\|\le \tau,$ the following recursive inequality:
\begin{equation}\label{rtq_0}\begin{aligned}
\|z_{n+1}-\theta_{0}\|^{2}&=\|z_{n}-\theta_{0}+z_{n+1}-z_{n}\|^{2}= \|z_{n}-\theta_{0}\|^{2}+2(z_{n}-\theta_{0})^{\top}(z_{n+1}-z_{n})+\|z_{n+1}-z_{n}\|^{2}.
\end{aligned}\end{equation}
Due to the definition of $z_{n+1}-z_{n},$ we have
\begin{equation}\nonumber\begin{aligned}
z_{n+1}-z_{n}&=\frac{u^{\top}(X_{n+1}-X_{n})-\alpha u^{\top}(X_{n}-X_{n-1})}{1-\alpha}\\&=\frac{-u^{\top}v_{n}+\alpha u^{\top}v_{n-1}}{1-\alpha}\\&=-\frac{\epsilon_{n}u^{\top}G(X_{n},\xi_{n})}{1-\alpha}.
\end{aligned}\end{equation}
Substitute above equation into Eq. \eqref{rtq_0}, and take the mathematical expectation, noting $\Expect(G(X_{n},\xi_{n}))=\Expect(G(X_{n})),$ getting
\begin{equation}\label{iuy_0}\begin{aligned}
\Expect\|z_{n+1}-\theta_{0}\|^{2}=\Expect\|z_{n}-\theta_{0}\|^{2}-\frac{2\epsilon_{n}}{1-\alpha}\cdot\Expect\big((z_{n}-\theta_{0})^{\top}u^{\top}G(X_{n})\big)+\frac{\epsilon_{n}^{2}}{(1-\alpha)^{2}}\Expect\|u^{\top}G(X_{n},\xi_{n})\|^{2}.
\end{aligned}\end{equation} For $u^{\top}G(X_{n}),$ dur to Eq. \eqref{mSGD_13}, we get
\begin{equation}\nonumber\begin{aligned}
u^{\top}G(X_{n})&=\nabla g(u^{\top}{X_{n}})+ (u^{\top}G(X_{n})-\nabla g(u^{\top}{X_{n}}))\\&=\nabla g(z_{n})+\frac{\alpha}{1-\alpha}(\nabla g(u^{\top}X_{n})-\nabla g(z_{n}))+(u^{\top}G(X_{n})-\nabla g(u^{\top}{X_{n}}))+.
\end{aligned}\end{equation}
Then we get
\begin{equation}\nonumber\begin{aligned}
-\frac{2\epsilon_{n}}{1-\alpha}\cdot\Expect\big((z_{n}-\theta_{0})^{\top}u^{\top}G(X_{n})\big)\le -\frac{2\epsilon_{n}}{1-\alpha}\cdot\Expect\big((z_{n}-\theta_{0})^{\top}\nabla g(z_{n})\big)+\mathcal{O}(\epsilon_{n}^{2}).
\end{aligned}\end{equation} Substitute above inequity into Eq. \eqref{rtq_0}, acquiring
\begin{equation}\label{iuy_1}\begin{aligned}
\Expect\|z_{n+1}-\theta_{0}\|^{2}=\Expect\|z_{n}-\theta_{0}\|^{2}-\frac{2\epsilon_{n}}{1-\alpha}\cdot\Expect\big((z_{n}-\theta_{0})^{\top}\nabla g(z_{n})\big)+\mathcal{O}(\epsilon_{n}^{2}).
\end{aligned}\end{equation}
For any term \(k\) in the first \(T\) iterations \(1, 2, \ldots, T\), we set \(\theta_{0}\) in Eq. \eqref{rtq_0} to \(z_{T-k}\), obtaining $\exists\ l>0,\ l_{0}>0$ such that
\begin{equation}\nonumber\begin{aligned}
\sum_{t=T-k}^{T}\Expect((z_{t}-z_{T-k})^{\top}\nabla g(z_{n}))\le \frac{l}{\sqrt{m}}\big(\sqrt{T}-\sqrt{T-k}\big)+l_{0}\sqrt{m}\sum_{t=T-k}^{T}\frac{1}{\sqrt{t}}.
\end{aligned}\end{equation}
By convexity, we can lower bound $(z_{t}-z_{T-k})^{\top}\nabla g(z_{t})$ bt $g(z_{t})-g(z_{T-k}).$ Also, it is easy to get that 
\[\sum_{t=T-k}^{T}\frac{1}{\sqrt{t}}\le 2(\sqrt{T}-\sqrt{T-k-1}).\] Then we get
\begin{equation}\nonumber\begin{aligned}
\Expect\Bigg(\sum_{t=T-k}^{T}(g(z_{t})-g(z_{T-k}))\bigg)\le \left(\frac{l}{\sqrt{m}}+l_{0}\sqrt{m}\right)\big(\sqrt{T}-\sqrt{T-k-1}\big)\le \left(\frac{l}{\sqrt{m}}+l_{0}\sqrt{m}\right)\frac{k+1}{\sqrt{T}}.
\end{aligned}\end{equation}
Then we define \(S_{k}=\frac{1}{k+1}\sum_{t=T-k}^{T}\Expect(g(z_{t}))\) be the expected average value of the last $K+1$ iterates. The bound above implies that
\[-\Expect(g(z_{T-k}))\le -\Expect(S_{k})+\frac{l\sqrt{m}+\frac{l_0}{\sqrt{m}}}{\sqrt{T}}.\] By the definition of $S_{k}$ and the inequity above, we have
\begin{equation}\nonumber\begin{aligned}
k\Expect(S_{k-1})=(k+1)\Expect(S_{k})-\Expect(g(z_{T-k}))\le (k+1)\Expect(S_{k})-\Expect(S_{k})+\frac{l\sqrt{m}+\frac{l_0}{\sqrt{m}}}{\sqrt{T}},
\end{aligned}\end{equation}
and dividing by $k,$ implies
\[\Expect(S_{k-1})\le \Expect(S_{k})+\frac{l\sqrt{m}+\frac{l_0}{\sqrt{m}}}{k\sqrt{T}}.\]
Using the inequity repeatedly and by summing over $k=1,...,T-1,$ we have
\[\Expect(g(z_{T}))=\Expect(S_{0})\le \Expect(S_{T-1})+\frac{l\sqrt{m}+\frac{l_0}{\sqrt{m}}}{\sqrt{T}}\sum_{k=1}^{T-1}\frac{1}{k}.\] Using Eq. \eqref{iuy_1} with $k=T-1$ and $\theta_{0}=\theta^*,$ we can get
\[\Expect(S_{T-1})-g(\theta^*)\le \frac{l\sqrt{m}+\frac{l_0}{\sqrt{m}}}{\sqrt{T}}.\]
Finally, we get
\[\Expect(g(u^{\top}X_{T})-g(\theta^{*}))=\mathcal{O}\left(\left(l\sqrt{m}+\frac{l_0}{\sqrt{m}}\right)\frac{\ln T}{\sqrt{T}}\right).\]
\end{proof}

\subsection{Proof of Theorem \ref{dacnie}}
We define another event

	\begin{equation*}
		B_{n}=\{\|\nabla g(\overline{x}_{1})\|^{2}> a_{0},\ \|\nabla g(\overline{x}_{2})\|^{2}> a_{0}\,\cdots,\|\nabla g(\overline{x}_{n})\|^{2}> a_{0}\},
	\end{equation*}
and its characteristic function as $I_{n}^{(a_{0})}$. Then through Assumption \ref{ass_g1} and $\epsilon_{n}\ge\epsilon_{n+1}$ we get that

	\begin{equation}\nonumber\begin{aligned}
			& I_{n+1}^{(a_{0})}{g}(\overline{x}_{n+1})-I_{n}^{(a_{0})}{g}(\overline{x}_{n})= -\frac{1}{2}\sum_{k=i}^{n}\alpha^{n-k}\epsilon_{k}I^{(a_{0})}_{k}\|\nabla g(\overline{x}_{k})\|^{2}+\frac{\hat{\mu}_{0}\sigma_{0}^{2}}{2}\sum_{k=i}^{n}\alpha^{n-k}I^{(a_{0})}_{k}O(\epsilon^{2}_{k})+\overline{k}\alpha^{n}+\zeta_{n}\\
   &+\hat{L}\sum_{k=1}^{n}\alpha^{n-k}I^{(a_{0})}_{k}\epsilon_{k}\beta_{k},
	\end{aligned}\end{equation}	
 where $\overline{k}$, $\hat{L}$ and $\hat{\mu}_{0}$ are three constants which can not affect the result. Notice that 
 
 	\begin{equation}\nonumber\begin{aligned}
 		I^{(a_{0})}_{k}O(\epsilon^{2}_{k})&\le \frac{1}{a_{0}}I^{(a_{0})}_{k}\|\nabla g(\overline{x}_{k})\|^{2}O(\epsilon^{2}_{k})\\
   &=I^{(a_{0})}_{k}\|\nabla g(\overline{x}_{k})\|^{2}O(\epsilon^{2}_{k}).
 	\end{aligned}\end{equation}	
 Then we get
 
	\begin{equation}\nonumber\begin{aligned}
	&	I_{n+1}^{(a_{0})}{g}(\overline{x}_{n+1})-I_{n}^{(a_{0})}{g}(\overline{x}_{n})= -\frac{1}{2}\sum_{k=i}^{n}\alpha^{n-k}\big(\epsilon_{k}-O(\epsilon_{k}^{2})\big)\Expect\big(I^{(a_{0})}_{k}\|\nabla g(\overline{x}_{k})\|^{2}\big)+\zeta_{n}.
	\end{aligned}\end{equation}	
Due to  $\Expect(\zeta_{n})=0$, we make the mathematical expectation to obtain

	\begin{equation}\nonumber\begin{aligned}
		&\Expect\big(I_{n+1}^{(a_{0})}{g}(\overline{x}_{n+1})\big)-\Expect\big(I_{n}^{(a_{0})}{g}(\overline{x}_{n})\big)= -\frac{1}{2}\sum_{k=i}^{n}\alpha^{n-k}\big(\epsilon_{k}-O(\epsilon_{k}^{2})\big)I^{(a_{0})}_{k}\|\nabla g(\overline{x}_{k})\|^{2}.
	\end{aligned}\end{equation}	
We denote
\begin{equation}\begin{aligned}\nonumber
		\hat{F}_{n}^{(a_{0})}=\sum_{i=1}^{n}\bigg(\frac{1}{2-\alpha}\bigg)^{n-i}\Expect\big(I_{n}^{(a_{0})}g(\theta_{n})\big).
\end{aligned}\end{equation}
For convenience, we let 

	\begin{equation}\nonumber\begin{aligned}
			&\hat{G}_{n}^{(a_{0})}=\frac{2}{(1-\alpha)^{2}}\sum_{i=1}^{n}\bigg(\frac{1}{2-\alpha}\bigg)^{n-i}(\epsilon_{i}-O(\epsilon^{2}_{i}))\\
   &\cdot\Expect\big(I_{i}^{(a_{0})}\|\nabla g(\overline{x}_{n})\|^{2}\big).
	\end{aligned}\end{equation}
Then we get
\begin{equation}\nonumber\begin{aligned}
		\hat{F}_{n+1}^{(a_{0})}-\hat{F}_{n}^{(a_{0})}\le-\hat{G}_{n}^{(a_{0})}, 	\end{aligned}\end{equation} so there is
\begin{equation}\begin{aligned}\nonumber
		&\hat{F}_{n+1}^{(a_{0})}\le \hat{F}_{n}^{(a_{0})}\Bigg(1-\frac{\hat{G}_{n}^{(a_{0})}}{\hat{F}_{n}^{(a_{0})}}\Bigg)\le \hat{F}_{1}^{(a_{0})}\prod_{i=1}^{n}\Bigg(1-\frac{\hat{G}_{i}^{(a_{0})}}{\hat{F}_{i}^{(a_{0})}}\Bigg)\le \hat{q}_{0}\prod_{i=1}^{n}\Bigg(1-\frac{\hat{G}_{i}^{(a_{0})}}{\hat{F}_{i}^{(a_{0})}}\Bigg),    
\end{aligned}\end{equation} where $\hat{q}_{0}$ is a constant. We focus on $\frac{\hat{G}_{i}^{(a_{0})}}{\hat{F}_{i}^{(a_{0})}}$. Using $O'stolz\ theorem$ yields
\begin{equation}\nonumber\begin{aligned}
		\liminf_{i\rightarrow+\infty}\frac{\hat{G}_{i}^{(a_{0})}}{\epsilon_{i}\hat{F}_{i}^{(a_{0})}}&=\liminf_{i\rightarrow+\infty}\frac{2}{(1-\alpha)^{2}}\frac{\sum_{t=1}^{i}(\frac{1}{2-\alpha})^{i-t}\Expect(I_{t}^{(a_{0})}\|\nabla g(\overline{x}_{t})\|^{2})}{\sum_{t=1}^{i}(\frac{1}{2-\alpha})^{i-t}\Expect(I_{t}^{(a_{0})}g(\overline{x}_{t}))}\\&\ge \liminf_{i\rightarrow+\infty}\frac{2}{(1-\alpha)^{2}}\frac{\Expect(I_{i}^{(a_{0})}\|\nabla g(\overline{x}_{i})\|^{2})}{\Expect(I_{t}^{(a_{0})}g(\overline{x}_{i}))}.
\end{aligned}\end{equation}In the setting of this theorem, the loss function is bounded. We let $g(x)<\hat{T}$. Then there is

	\begin{equation}\nonumber\begin{aligned}
		\liminf_{i\rightarrow+\infty}\frac{2}{(1-\alpha)^{2}}\frac{\Expect(I_{i}^{(a_{0})}\|\nabla g(\overline{x}_{i})\|^{2})}{\Expect(I_{t}^{(a_{0})}g(\overline{x}_{i}))}\ge \liminf_{i\rightarrow+\infty}\frac{2}{(1-\alpha)^{2}}\frac{a}{\hat{T}}.
	\end{aligned}\end{equation}

Then it holds that  
\begin{equation}\begin{aligned}\nonumber
		\Expect(I_{n}^{(a_{0})})\le\hat{F}_{n+1}^{(a_{0})}=O\Big(e^{-\frac{s}{(1-\alpha)^{2}}\sum_{i=1}^{n}\epsilon_{n}}\Big),   
\end{aligned}\end{equation}
where $s=2a/\hat{T}$.

\end{document}